\documentclass[10pt,a4paper]{article}
\usepackage[british]{babel}
\usepackage{amssymb,amsmath,amsthm}
\usepackage{cite}
\usepackage{graphicx}
\usepackage[bookmarks=true, bookmarksopen=true, bookmarksopenlevel=4]{hyperref}
\usepackage{color}
\usepackage{enumerate}
\usepackage{pgf,tikz}
\usepackage{soul}
\usepackage{float}
\usepackage{subfig}
\usetikzlibrary{decorations.pathreplacing,matrix,calc,arrows,shapes.geometric}

\hypersetup{pdftitle={Spurious Quasi-Resonances},
  pdfauthor={R. Hiptmair, A. Moiola, E.A. Spence},
  colorlinks=true, linkcolor=myblue,  citecolor=myblue, filecolor=myblue,   urlcolor=myblue}

\allowdisplaybreaks[4]
\definecolor{myblue}{rgb}{0,0,0.6} 
\newcommand*{\der}[2]{\frac{\partial #1}{\partial #2}}
\newcommand{\iin}{\;\text{in}\;}
\newcommand{\oon}{\;\text{on}\;}

\newcommand*{\N}[1]{\left\|#1\right\|}

\newcommand{\Uu}[1]{{\mathbf{#1}}}
\newcommand{\IR}{\mathbb{R}}
\newcommand{\bn}{{\Uu n}}\newcommand{\bx}{{\Uu x}}\newcommand{\by}{{\Uu y}}

\newtheorem{theorem}{Theorem}[section]
\newtheorem{defin}[theorem]{Definition}
\newtheorem{lemma}[theorem]{Lemma}
\newtheorem{cor}[theorem]{Corollary}

\theoremstyle{definition}
\newtheorem{remark}[theorem]{Remark}
\newtheorem{example}[theorem]{Example}
\newcommand{\ee}{{\rm e}}
\newcommand{\ri}{{\rm i}}
\newcommand{\rd}{{\rm d}}
\newcommand{\beq}{\begin{equation}}      \newcommand{\eeq}{\end{equation}}
\newcommand{\beqs}{\begin{equation*}}    \newcommand{\eeqs}{\end{equation*}}
\newcommand{\bit}{\begin{itemize}}       \newcommand{\eit}{\end{itemize}}
\newcommand{\ben}{\begin{enumerate}}     \newcommand{\een}{\end{enumerate}}
\newcommand{\bal}{\begin{align}}         \newcommand{\eal}{\end{align}}
\newcommand{\bals}{\begin{align*}}       \newcommand{\eals}{\end{align*}}
\newcommand{\bse}{\begin{subequations}}	 \newcommand{\ese}{\end{subequations}}
\newcommand{\bpr}{\begin{proposition}}   \newcommand{\epr}{\end{proposition}}
\newcommand{\bre}{\begin{remark}}        \newcommand{\ere}{\end{remark}}
\newcommand{\bpf}{\begin{proof}}         \newcommand{\epf}{\end{proof}}
\newcommand{\ble}{\begin{lemma}}         \newcommand{\ele}{\end{lemma}}
\newcommand{\bco}{\begin{corollary}}     \newcommand{\eco}{\end{corollary}}
\newcommand{\bex}{\begin{example}}       \newcommand{\eex}{\end{example}}
\newcommand{\bth}{\begin{theorem}}       \newcommand{\enth}{\end{theorem}}
\newcommand{\Rea}{\mathbb{R}}            \newcommand{\Com}{\mathbb{C}}
\newcommand{\Oi}{\Om}               \newcommand{\Oe}{\Op}
\newcommand{\Om}{\Omega^-}               \newcommand{\Op}{\Omega^+}

\newcommand{\pdiff}[2]{\frac{\partial #1}{\partial #2}}

\newcommand{\supp}{\operatorname{supp}}

\newcommand{\half}{\frac{1}{2}}

\newcommand{\tendi}{\rightarrow \infty}

\newcommand{\Oin}{{\Om}}
\newcommand{\Oout}{{\Op}}
\newcommand{\nin}{{n_i}}
\newcommand{\nout}{{n_o}}

\newcommand{\SRC}{\mathrm{SRC}}
\newcommand{\loc}{_{\mathrm{loc}}}

\newcommand{\mythmname}[1]{\textbf{\emph{(#1.)}}}

\newcommand{\cH}{{\mathcal H}}
\newcommand{\cN}{{\mathcal N}}

\newcommand{\cK}{{\mathcal K}}
\newcommand{\cV}{{\mathcal V}}

\newcommand{\tin}{\text{ in }}
\newcommand{\tfa}{\text{ for all }}
\newcommand{\tfor}{\text{ for }}
\newcommand{\tas}{\text{ as }}
\newcommand{\tand}{\text{ and }}

\numberwithin{equation}{section}

\usepackage{fancyhdr}
\newcommand{\Cauchy}{\boldsymbol{\gamma}_C}
\newcommand{\bphi}{\boldsymbol{\phi}}
\newcommand{\bphit}{\widetilde{\boldsymbol{\phi}}}
\newcommand{\bpsi}{\boldsymbol{\psi}}
\newcommand{\First}{A_{\rm I}}
\newcommand{\Second}{A_{\rm II}}
\newcommand{\Firstgen}{A^{\rm gen}_{\rm I}}
\newcommand{\AugFirst}{\widetilde{A}_{\rm I}}
\newcommand{\AugSecond}{\widetilde{A}_{\rm II}}
\newcommand{\Helmspace}{\mathsf{H}}

\definecolor{amcol}{rgb}{0.8,0,0}

\definecolor{escol}{rgb}{0,0,0.8}
\definecolor{estcol}{rgb}{0,0.8,0}

\title{Spurious Quasi-Resonances in Boundary Integral Equations for the Helmholtz Transmission Problem}

\author{
  Ralf Hiptmair\thanks{Seminar for Applied Mathematics, ETH Z\"urich, Z\"urich, Switzerland (\texttt{hiptmair@sam.math.ethz.ch})}, \,
  Andrea Moiola\thanks{
    Department of Mathematics, University of Pavia, Pavia, Italy (\texttt{andrea.moiola@unipv.it})},
  \, Euan A.~Spence\thanks{Department of Mathematical Sciences, University of Bath, Bath, UK (\texttt{E.A.Spence@bath.ac.uk})}}

\date{\today}

\begin{document}
\maketitle

\begin{abstract}
We consider the Helmholtz transmission problem 
with piecewise-constant material coefficients, and the standard associated direct boundary integral equations.
For certain coefficients and geometries, the norms of the inverses of the boundary integral operators grow rapidly through an increasing sequence of frequencies, even though this is not the case for the solution operator of the transmission problem; we call this phenomenon that of  \emph{spurious quasi-resonances}. We give a rigorous explanation of why and when spurious quasi-resonances occur, and propose modified boundary integral equations that are not affected by them. 


\end{abstract}

\section{Introduction and statement of the main results}
\label{sec:intro}

The goal of this paper is to explain, and also provide a remedy for, the feature of \emph{spurious quasi-resonances} in boundary integral equations for the Helmholtz transmission problem. This feature is illustrated in numerical experiments in \S\ref{sec:spqr}, with the explanation and our remedy given  in \S\ref{sec:mainresults}. \S\S\ref{sec:transmission}-\ref{sec:Calderon} define, respectively, the Helmholtz transmission problem, its solution operator, and the standard first- and second-kind direct boundary integral formulations of this transmission problem.

\subsection{The Helmholtz transmission scattering problem}
\label{sec:transmission}

We consider the scattering of an incident time-harmonic acoustic wave by a penetrable
  homogeneous object that occupies the region of space $\Oin\subset\IR^d$, $d=2,3$, which
  is a bounded Lipschitz open set. We first introduce notation
  necessary for a precise mathematical statement of this \emph{transmission problem}. Let $\Oout:=\IR^d\setminus\overline\Oin$, $\Gamma:=\partial \Om=\partial \Op$, and let $\bn$ be the unit normal vector field on $\Gamma$ pointing from $\Oin$ into $\Oout$. For any
$\varphi\in L^2\loc(\IR^d)$, we let $\varphi^-:=\varphi|_\Oin$ and
$\varphi^+:=\varphi|_\Oout$. 
With $H^1_{\rm loc}(\Omega^\pm,\Delta):=\{ v : \chi v  \in H^1(\Omega^\pm), \Delta (\chi v) \in L^2(\Omega^\pm) \tfa \chi \in C^\infty_{\rm comp}(\Rea^d)\}$, we define the Dirichlet and Neumann trace operators
\begin{equation*}
  \gamma_D^{\pm} : H^1_{\rm loc}(\Omega_{\pm}) \rightarrow H^{1/2}(\Gamma) \quad\tand\quad
  \gamma_N^{\pm} : H^1_{\rm loc}(\Omega_{\pm}, \Delta) \rightarrow H^{-1/2}(\Gamma),
\end{equation*}
with $\gamma^\pm_{D}v := v^\pm|_{\Gamma}$ and $\gamma_N^{\pm}$ such that if
$v\in H^2_{\rm loc}(\Omega_{\pm})$ then
$\gamma_N^{\pm} v = \bn \cdot \gamma_D^{\pm}(\nabla v)$.  Let
$\Cauchy^\pm := (\gamma_D^{\pm}, \gamma_N^\pm)$ be the Cauchy trace, which
satisfies
\begin{equation*}
  \Cauchy^{\pm}: H^1_{\rm loc}(\Omega^{\pm},\Delta) \rightarrow
  H^{1/2}(\Gamma)\times H^{-1/2}(\Gamma).
\end{equation*}
Given $\varphi\in C^1(\IR^d\setminus B_R)$, for some ball $B_R:=\{|\bx|<R\}$,  and $\kappa>0$, 
$\varphi$ satisfies the \emph{Sommerfeld radiation condition} if 
\begin{equation}
  \label{eq:src}
  \lim_{r\to\infty}r^{\frac{d-1}2}\left(\der{\varphi(\bx)}r-\ri\kappa \varphi(\bx)\right)=0
\end{equation}
uniformly in all directions, where $r:=|\bx|$; we then write $\varphi\in\SRC(\kappa)$.

Given $\nin,\nout>0$ and frequency $k>0$, the \textbf{Helmholtz transmission scattering problem} is that of finding the complex amplitude $u$ of the sound
  pressure, with $u\in H^1\loc(\IR^d\setminus\Gamma)$ the solution of
\begin{align}
  \begin{aligned}
    (\Delta +k^2 \nin)u^- &=0&&\iin\Oin,\\
    (\Delta +k^2 \nout)u^+ &=0 &&\iin\Oout,\\
    \Cauchy^- u^- &= \Cauchy^+ u^+ + \Cauchy^{\pm} u^I 
    &&\oon \Gamma,\\
    u^+ &\in\SRC(k\sqrt{n_o}),
  \end{aligned}
  \label{eq:htsp}
\end{align}  
where the incident wave $u^{I}$ is an entire solution of the homogeneous Helmholtz
equation in $\Rea^{d}$, 
\begin{equation}
  \label{eq:entire}
  (\Delta + k^2 \nout) u^I=0 \quad\tin \Rea^d.
\end{equation}
This set up means that $u^-$ is the total field in $\Oi$ and $u^+$ the scattered field in $\Oe$.

In principle, the jump $\Cauchy^+ u^+-\Cauchy^- u^-$ of the Cauchy trace of $u$ across
$\Gamma$ can be more general than the Cauchy trace of an incident wave. This leads to the
following generic  Helmholtz transmission problem.

\begin{defin}\mythmname{The Helmholtz transmission problem}\label{def:HTP}
  Given 
  positive real numbers $k, \nin,$ and $\nout$ and ${\bf f} \in H^{1/2}(\Gamma)\times H^{-1/2}(\Gamma)$, 
  find $u\in H^1\loc(\IR^d\setminus\Gamma)\cap\SRC(k\sqrt{n_o})$ such that,
  \begin{align}
    \begin{aligned}
      (\Delta +k^2 \nin)u^- &=0&&\iin\Oin,\\
      (\Delta +k^2 \nout)u^+ &=0 &&\iin\Oout,\\
      \Cauchy^- u^- &= \Cauchy^+ u^+ +
      {\bf f}
      &&\oon \Gamma.
    \end{aligned}
    \label{eq:BVP}
  \end{align}
\end{defin}

The following well-posedness result is proved in, e.g., \cite[Lemma 2.2 and Appendix A]{MoSp:19}.

\ble\label{lem:wellposed}
The solution of the transmission problem of Definition~\ref{def:HTP} exists and is unique. Moreover, 
if ${\bf f}\in H^1(\Gamma)\times L^2(\Gamma)$ then
$\Cauchy^{\pm}u^{\pm} \in H^1(\Gamma)\times L^2(\Gamma)$.
\ele

\begin{remark}
  \label{rem:genform}
  The transmission problem of Definition \ref{def:HTP} is not the most general form of the
  transmission problem. If the transmission condition in \eqref{eq:BVP} is replaced by
  \beq\label{eq:newtc} \Cauchy^- u^- = D \Cauchy^+ u^+ +{\bf f}, \quad \text{ where } D:=
  \left(
    \begin{array}{cc}
      1 & 0 \\
      0 & \alpha 
    \end{array}
  \right) \eeq for $\alpha$ a constant, then this covers all possible constant-coefficient
  transmission problems; see, e.g., \cite[Page 322]{MoSp:19}. In Appendix
  \ref{app:general} we outline how our results extend this more general transmission
  problem. We see that, although the main ideas remain the same, more notation and
  technicalities are required, hence why we have chosen to focus on the simpler problem of
  Definition \ref{def:HTP}.
\end{remark}

\subsection{Solution operators and quasi-resonances}\label{sec:Sio}

\begin{defin}\mythmname{Solution operators}\label{def:Sio}
  Given positive real numbers $k, c_i,$ and $c_o$, let
  \begin{equation*}
    S(c_i,c_o){\bf f}:= \Cauchy^- u,
  \end{equation*}
  where $u\in H^1_{\rm loc}(\Rea^d\setminus\Gamma)\cap \SRC(k\sqrt{c_o})$ is the solution of 
the Helmholtz transmission problem
  \begin{align}
    \begin{aligned}
      (\Delta +k^2 c_i)u^- &=0&&\iin\Oin,\\
      (\Delta +k^2 c_o)u^+ &=0 &&\iin\Oout,\\
      \Cauchy^- u^- &= \Cauchy^+ u^+ +{\bf f} &&\oon \Gamma.
    \end{aligned}
    \label{eq:BVP2}
  \end{align}
\end{defin}

Lemma \ref{lem:wellposed} implies that $S(c_i,c_o)$ is well defined and bounded on
\emph{either} $H^{1/2}(\Gamma) \times H^{-1/2}(\Gamma)$ \emph{or}
$H^1(\Gamma)\times L^2(\Gamma)$.  We introduce the abbreviations
\begin{equation*}
  S_{io}:= S(n_i,n_o) \quad\tand\quad S_{oi}:= S(n_o,n_i).
\end{equation*}
We refer to $S_{io}$ as the ``physical'' solution operator, since it corresponds to the
transmission problem of Definition \ref{def:HTP}, and $S_{oi}$ as the ``unphysical''
solution operator, since it corresponds to the transmission problem where the indices
$n_i$ and $n_o$ are swapped compared to those in Definition~\ref{def:HTP}.
The results of this paper show that to understand the behaviour of boundary integral operators used to solve the ``physical'' problem, one needs \emph{both} the ``physical'' solution operator \emph{and} the ``unphysical'' one (this is made more precise in Theorem \ref{thm:1} below).

The spurious quasi-resonances we study in this paper are related to the high-frequency behaviour of boundary integral operators. 
We therefore recap here the high-frequency behaviour of $S_{io}$; recall that this 
 depends on which of $n_i$ and $n_o$
is larger.  Indeed, if $n_i<n_o$ and $\Oi$ is Lipschitz and star-shaped with respect to a
ball, then Lemma \ref{lem:bound2} below shows that the norm of $S_{io}$ has, at worst,
mild algebraic growth in $k$; this result uses the bounds on the solution operator from
\cite{MoSp:19}, with analogous bounds obtained for smooth, convex $\Oi$ with strictly
positive curvature in \cite{CaPoVo:99}. If $n_i> n_o$ and $\Oi$ is smooth and convex with
strictly positive curvature, then Lemma \ref{lem:bound3} below, based on the results of
\cite{PoVo:99}, shows that there exists $0<k_1<k_2<\ldots$ with $k_j\tendi$ such that the
norm of $S_{io}$ blows up superalgebraically through $k_j$ as $j\tendi$. (Similar results
in the particular case when $\Oi$ is a ball were obtained in \cite{Cap12, CLP12}, and
summarised in \cite[Chapter~5]{AC16}).

We call these real frequencies $k_j$ \emph{quasi-resonances}, since they can be understood
as real parts of complex resonances of the transmission problem lying close to the real
axis (with this terminology also used in, e.g., \cite{Cap12, CLP12, AC16}); the particular
functions on which $S_{io}$ at $k=k_j$ blows up are called \emph{quasimodes}. The
relationship between quasimodes and resonances is a classic topic in scattering theory;
see \cite{TaZw:98, St:99, St:00}, \cite[\S7.3]{DyZw:19}. 
The Weyl-type bound on the number of resonances of the transmission problem when $\Oi$ is
smooth and convex with strictly positive curvature in \cite[Theorem 1.3]{CaPoVo:01}
implies that the number of quasi-resonances in $[0,K]$ in this case grows like $K^d$ as
$K\tendi$.

\begin{remark}
  \label{rem:qr}
  The physical reason for the existence of quasi-resonances when $n_i>n_o$ is that, in
  this case, geometric-optic rays can undergo total internal reflection when hitting
  $\Gamma$ from $\Omega^-$. Rays ``hugging" the boundary via a large number of bounces
  with total internal reflection correspond to solutions of the transmission problem
  localised near $\Gamma$; in the asymptotic-analysis literature these solutions are known
  as ``whispering gallery" modes; see, e.g., \cite{BaBu:91, BaDaMo:21}.
  The existence of quasi-resonances of the transmission problem has only been rigorously
  established when $\Oi$ is smooth and convex with strictly positive curvature. The
  understanding above via rays suggests that such quasi-resonances and quasimodes do not
  exist for polyhedral $\Oi$ (since sharp corners prevent rays from moving parallel to the
  boundary), although solutions with localisation qualitatively similar to that of
  quasimodes can be seen when $\Oi$ is a pentagon \cite[Figure 13]{LeDjArLaZyDuScBo:07} or
  a hexagon \cite[Figure 23]{CaWi:15}.
\end{remark}

\subsection{Calder\'on projectors and the standard first- and second-kind direct boundary integral equations (BIEs)}\label{sec:Calderon}

Since all the layer potentials and integral operators depend on $k$, we omit this $k$-dependence in the notation. 
Let the Helmholtz fundamental solutions be given by 
\begin{align*}
  \Phi_{i/o} (\bx,\by) := 
  \dfrac{\ri}{4} H^{(1)}_0(k\sqrt{n_{i/o}} \lvert \bx - \by \rvert )\quad \tfor d=2,\,\, \tand\,\,
  \Phi_{i/o} (\bx,\by) :=             \dfrac{\ee^{\ri k \sqrt{n_{i/o}}\lvert \bx-\by \rvert}}{4 \pi \lvert \bx - \by \rvert}\quad \tfor d=3,
\end{align*}
where \(H^{(1)}_0\) is the Hankel function of the first kind and order zero; see  \cite[Section~3.1]{SauterSchwab}.

As in \cite[Equation~3.6]{SauterSchwab}, the single-layer, adjoint-double-layer,
double-layer, and hypersingular operators are defined for $\phi\in L^2(\Gamma)$ and
$\psi\in H^1(\Gamma)$ by
\begin{align}\label{eq:SD'}
  &V_{i/o} \phi(\bx) := \int_\Gamma \Phi_{i/o}(\bx,\by) \phi(\by)\,\rd s(\by), \qquad
  K'_{i/o} \phi(\bx) := \int_\Gamma \frac{\partial \Phi_{i/o}(\bx,\by)}{\partial n(\bx)}  \phi(\by)\,\rd s(\by),\\
  &K_{i/o} \phi(\bx) := \int_\Gamma \frac{\partial \Phi_{i/o}(\bx,\by)}{\partial n(\by)}  \phi(\by)\,\rd s(\by), 
  \quad 
  W_{i/o} \psi(\bx) := -\pdiff{}{n(\bx)} \int_\Gamma \frac{\partial \Phi_{i/o}(\bx,\by)}{\partial n(\by)}  \psi(\by)\,\rd s(\by),
  \label{eq:DH}
\end{align}
for $\bx \in \Gamma$ (note that the sign of the hypersingular operator is swapped compared
to, e.g., \cite{CGLS12}). When $\Gamma$ is Lipschitz, the integrals defining $K_{i/o}$ and
$K_{i/o}'$ must be understood as Cauchy principal values (see, e.g., \cite[Equation
2.33]{CGLS12}), and the integral defining $W_{i/o}$ is understood as a non-tangential
limit (see, e.g., \cite[Equation 2.36]{CGLS12}) or finite-part integral (see, e.g., \cite[Theorem 7.4 (iii)]{MCL00}), but we do not
need the details of these definitions in this paper.

Let the Calder\'on projectors $P^\pm_{i/o}$ be defined by 
\begin{equation}\label{eq:Pdef}
  P^{\pm}_{i/o} := \half I \pm M_{i/o},
\quad\text{ where }\quad
  M_{i/o}:= 
  \left[
    \begin{array}{cc}
      K_{i/o} & - V_{i/o}\\
      -W_{i/o}& -K_{i/o}'
    \end{array}
  \right];
\end{equation}
see, e.g., \cite[Section~3.6]{SauterSchwab}, \cite[Page 117]{CGLS12}.
Basic results about $P^\pm_{i/o}$ (including that they are indeed projectors) are in
\S\ref{sec:background}, but we record here that
\begin{equation}\label{eq:sum}
  P^+_{i/o}+ P^-_{i/o} = I.
\end{equation}

Let the boundary integral operators (BIOs) $\First$ and $\Second$ be defined by
\begin{align}\label{eq:First}
  \First 
  &:= P^-_o-P^+_i = P^-_i - P^+_o=
  \left[
    \begin{array}{cc}
      -(K_i+ K_o) & V_i+ V_o\\
      W_i+ W_o & K'_i+ K'_o
    \end{array}
  \right]
\end{align}
and
\begin{align}\label{eq:Second}
  \Second 
  &:= P^-_o + P^+_i = 2I -P^+_o-P^-_i
  =
  I +
  \left[
    \begin{array}{cc}
      K_i- K_o & -(V_i- V_o)\\
      -(W_i- W_o) & -(K'_i- K'_o)
    \end{array}
  \right].
\end{align}

\ble\label{lem:3}
If $u$ is the solution of the Helmholtz transmission problem of Definition \ref{def:HTP}, then 
\begin{equation}\label{eq:firstsecondkind}
  \First (\Cauchy^- u^-) = P_o^-{\bf f}
  \quad\tand\quad
  \Second(\Cauchy^- u^-) =  P_o^-{\bf f}.
\end{equation}
In particular, if $u$ solves the Helmholtz transmission scattering problem
\eqref{eq:htsp}, then
\begin{equation}\label{eq:firstsecondkind2}
  \First (\Cauchy^- u^-) = \Cauchy^- u^I,
  \quad\tand\quad
  \Second(\Cauchy^- u^-) =  \Cauchy^- u^I.
\end{equation}
\ele

These boundary integral equations (BIEs) are called single-trace formulations (STFs).
  The \emph{first-kind} BIEs in \eqref{eq:firstsecondkind} and \eqref{eq:firstsecondkind2}
appeared in \cite{CoSt:85}, \cite{vo:89}, and
  are also derived in, e.g., \cite[Section~3.3]{ClHiJePi:15}. Their counterparts for electromagnetic
  scattering are known as the PMCHWT (Poggio--Miller--Chang--Harrington--Wu--Tsai)
  formulation \cite{POM73}. The \emph{second-kind} BIEs in \eqref{eq:firstsecondkind} and
  \eqref{eq:firstsecondkind2} can be found in, e.g., \cite{CHS12} and are known as the M\"uller
  formulation in computational electromagnetics \cite{MUE69}.

\ble\label{lem:invertible}
(i) Both $\First$ and $\Second$ are bounded and invertible on $H^{1/2}(\Gamma)\times H^{-1/2}(\Gamma)$. 

(ii) $\Second$ is bounded and invertible on $H^{1}(\Gamma)\times L^2(\Gamma)$.
\ele

The proofs of Lemmas \ref{lem:3} and \ref{lem:invertible} are contained in \S\ref{sec:background}.

The reason for the choice of spaces in Lemma \ref{lem:invertible} is the following. From the point of view of computation, the natural space in which to consider $\First$ is the trace space $H^{1/2}(\Gamma)\times H^{-1/2}(\Gamma)$, and the natural space in which to consider $\Second$ is the $L^2$-based space $H^{1}(\Gamma)\times L^2(\Gamma)$ (see, e.g., the discussion in \cite{CHS12} and the references therein); these choices are both included in Lemma \ref{lem:invertible}. 
It turns out that all the results for $\Second$ on $H^{1}(\Gamma)\times L^2(\Gamma)$ also hold on $H^{1/2}(\Gamma)\times H^{-1/2}(\Gamma)$, and thus we include this second choice of space for $\Second$.

\subsection{Spurious quasi-resonances for the standard BIOs}
\label{sec:spqr}

Lemma \ref{lem:invertible} shows that the BIEs of \eqref{eq:firstsecondkind} and \eqref{eq:firstsecondkind2} are well-posed. It is then reasonable to believe that the solution operators of these BIEs inherit the behaviour (with respect to frequency) of the solution operator of the transmission problem.
The following numerical results, however, show that this is not the case.
\footnote{The code used to produce the numerical results is available at \url{https://github.com/moiola/TransmissionBIE-OpNorms}}

\begin{example}
  \label{ex:numexp}
  If $\Gamma$ is a circle for $d=2$ or a sphere for $d=3$ all boundary integral operators
  $V_{i/o}$, $K'_{i/o}$, $K_{i,o}$, and $W_{i/o}$ can be ``diagonalized'' by switching to
  a ``modal'' $L^{2}(\Gamma)$-orthogonal basis of Fourier harmonics in 2D or spherical
  harmonics in 3D, respectively. The corresponding eigenvalues can be found in, e.g., 
  \cite{AMM98} for $d=2$ and in, e.g., \cite{VGG14} for $d=3$. All relevant norms have a simple
  sum representation with respect to these bases.  Therefore we can compute the norms of the
  solution operators as the maximum of the Euclidean norms of $2\times 2$-matrices, one
  for every mode. We did this in MATLAB for the modes of order at most 100, which seems to be
  sufficient, because the maximal norm was invariably found among the modes of order $\leq 25$.

  We report the computed norms of the solution operator $S_{io}$ along with the norms of $\First^{-1}$ and $\Second^{-1}$ (i.e., the solution operators for the BIEs
  \eqref{eq:firstsecondkind2}) on the space $H^{1/2}(\Gamma)\times H^{-1/2}(\Gamma)$, where we use the weighted norm $\|\cdot\|_{H^{1/2}_k(\Gamma)\times H^{-1/2}_k(\Gamma)}$ defined in \S\ref{sec:weighted}. We plot these norms for different frequencies $k$ and give the
  results for $d=2$ in Figure~\ref{fig:2d} and for $d=3$ in
  Figure~\ref{fig:3d}.

  \begin{figure}[H]
\includegraphics[width=0.5\textwidth,clip]{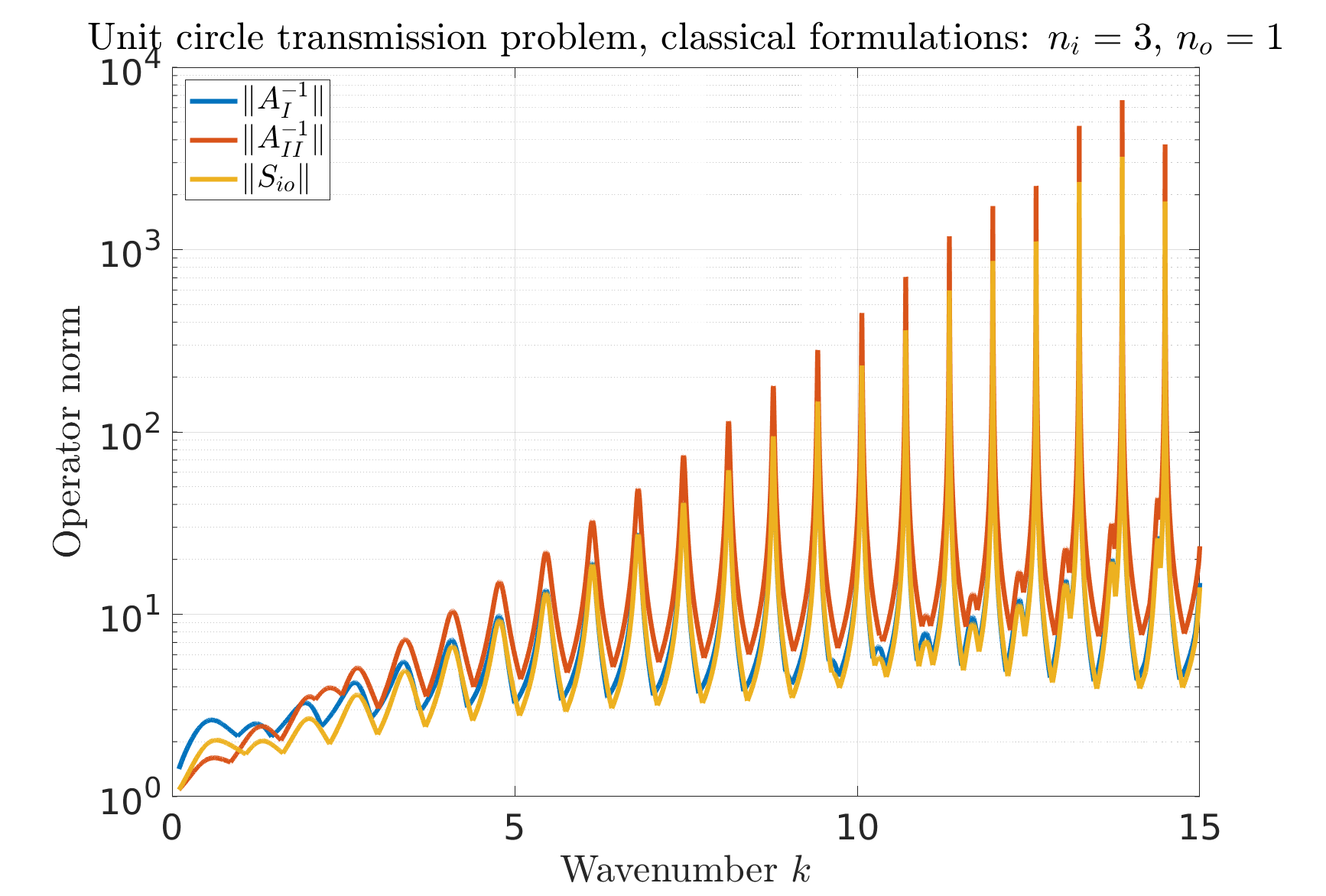}\hfill
\includegraphics[width=0.5\textwidth,clip]{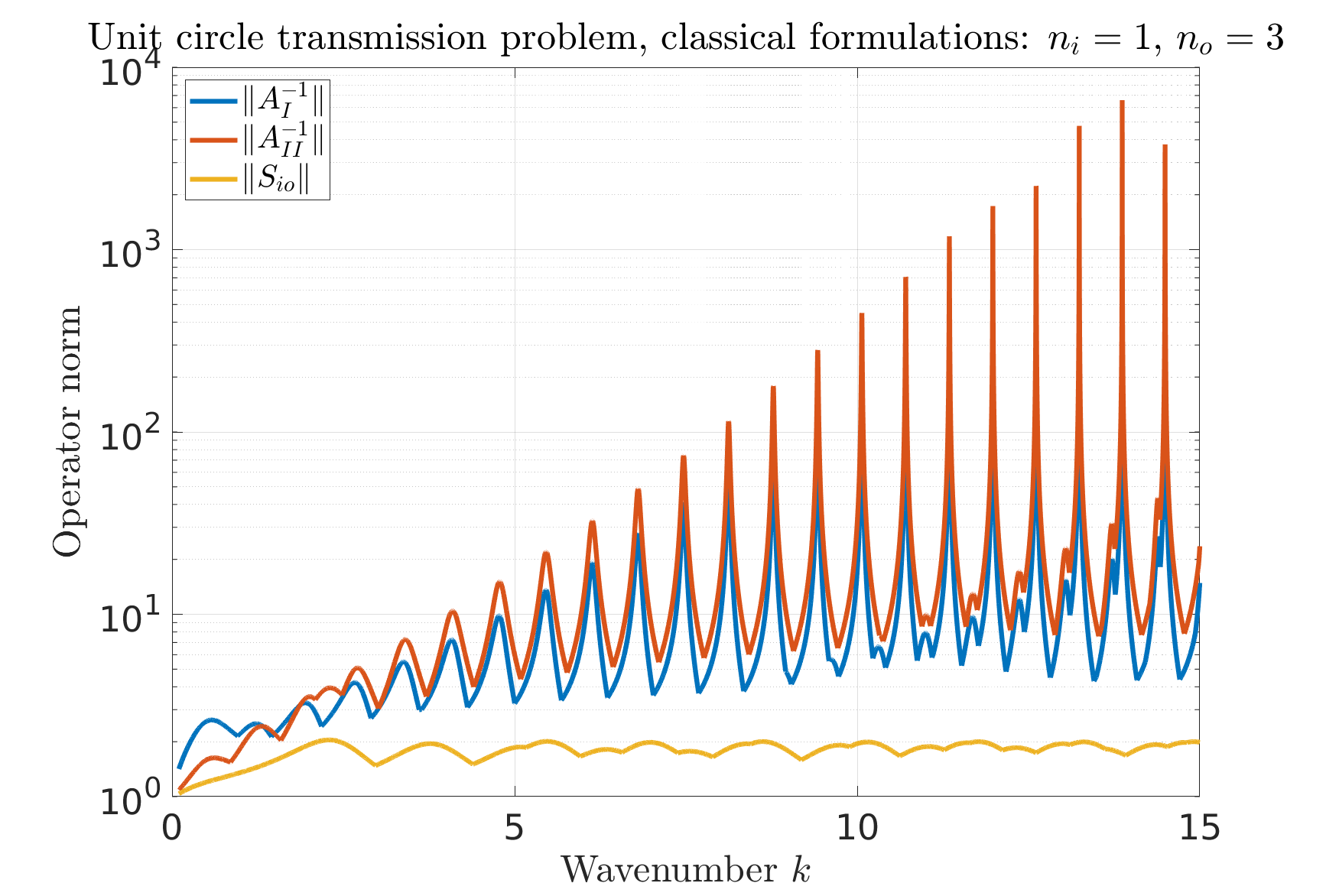}
    \caption{$\Gamma$ unit circle: Norms of operators $S_{io}$, $\First^{-1}$, and $\Second^{-1}$ on
      $H^{1/2}(\Gamma)\times H^{-1/2}(\Gamma)$ for $n_{i}=3$, $n_{o}=1$ (left) and
      $n_{i}=1$, $n_{o}=3$ (right)}
    \label{fig:2d}
  \end{figure}
  
   \begin{figure}[H]
\includegraphics[width=0.5\textwidth,clip]{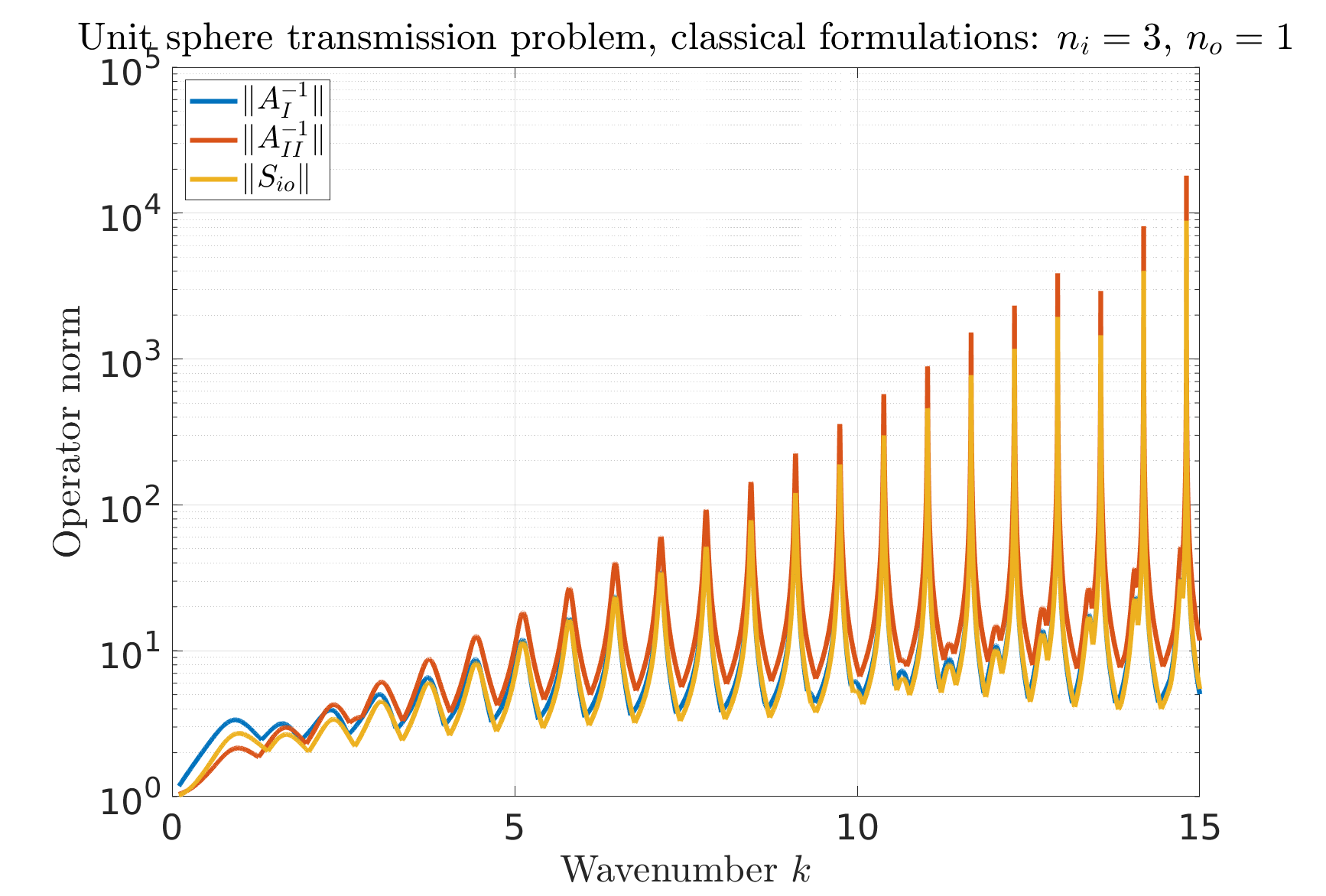}\hfill
\includegraphics[width=0.5\textwidth,clip]{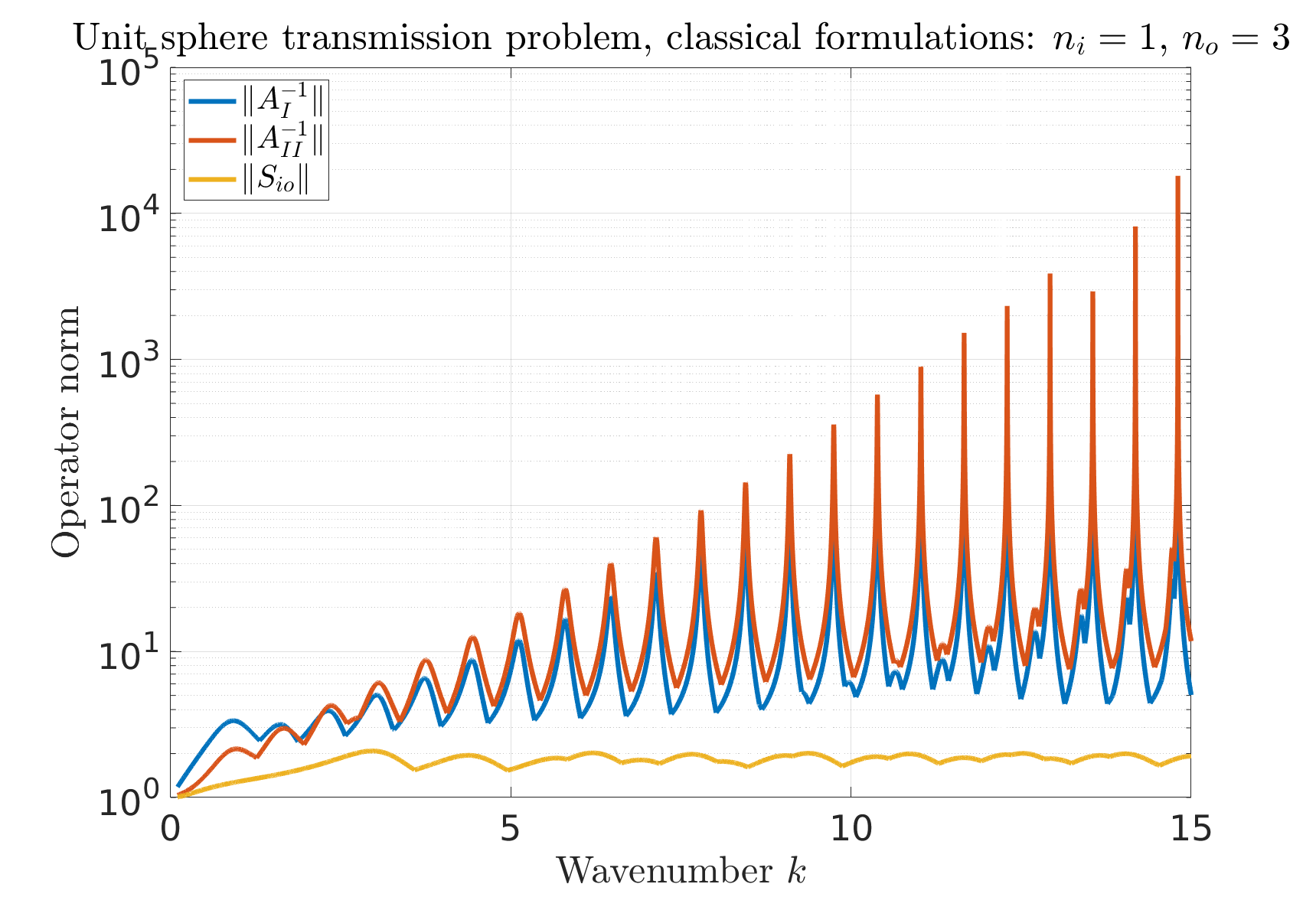}
    \caption{$\Gamma$ unit sphere: Norms of operators $S_{io}$, $\First^{-1}$, and $\Second^{-1}$ on
      $H^{1/2}(\Gamma)\times H^{-1/2}(\Gamma)$ for $n_{i}=3$, $n_{o}=1$ (left) and
      $n_{i}=1$, $n_{o}=3$ (right)}
    \label{fig:3d}
  \end{figure}

When $n_{i}>n_{o}$ (plots on the left) we see the typical spikes in the norms as a function of $k$, expected because of the results recalled in \S\ref{sec:Sio}. Indeed, the results in \S\ref{sec:Sio} predict super-algebraic growth through quasi-resonances only for sufficiently-large $k_j$. However, noting the logarithmic scale on the $y$-axis of the plots, we see that the super-algebraic growth occurs through the spikes even for small- to moderate-sized $k_j$.

    Conversely, for $n_{i}<n_{o}$
  (right plots) the norm of $S_{io}$ (in yellow) does not have any spikes, whereas the
  spikes persist in the norms of $\First^{-1}$ and $\Second^{-1}$. 
\end{example}

The observations made in Example~\ref{ex:numexp} provide evidence of \emph{spurious
  quasi-resonances} of $\First$ and $\Second$ when $n_{i}<n_{o}$: for certain frequencies these boundary integral operators are ill-conditioned though for the same frequencies the solution operator is stable. 
  
On rare occasions such spurious quasi-resonances have been noticed before. Indeed, the paper \cite{MiNiNi:17} computed the complex eigenvalues of $\First$ and $\Second$ and pointed out in \cite[Section~2.3]{MiNiNi:17} the existence of ``fictitious eigenvalues'', i.e., 
non-physical poles of the resolvent operators. Although \cite{MiNiNi:17} did not give a rigorous explanation for this
phenomenon,  \cite{MiNiNi:17} attempted to remedy it by modifying the BIEs; these new BIEs, however, still have issues with poles with small imaginary part -- see the discussion in \cite[\S4]{MiNiNi:17}. Non-physical spikes in the condition numbers of discretized BIEs for Helmholtz transmission problems were also
reported in \cite[Section~4.4]{WHG21}, but no deeper investigation was attempted.

The observation of the spurious quasi-resonances of Example~\ref{ex:numexp} was the starting point for this
paper -- we wanted to understand precisely why they affect $\First$ and $\Second$. We
also wanted to find alternative BIEs immune to spurious quasi-resonances. The remainder of this paper reports our progress towards these goals. 

\begin{remark}
  \label{rem:cfie}
  For the standard first and second-kind BIEs for the exterior Dirichlet and Neumann problems
  for the Helmholtz operator (modelling acoustic scattering by impenetrable objects), the
  occurrence of spurious (true) resonances is well-known; see, e.g.,  \cite[Section~3.9.2]{SauterSchwab}: the solutions of the BIEs are not unique for an infinite
  sequence of distinct $k$s, although the boundary-value problems have
  unique solutions for all $k$. 
The standard remedies for this 
are recalled (and linked to the results of the present paper) in Remark \ref{rem:impenetrable} below.
\end{remark}

\subsection{Statement of the main results}\label{sec:mainresults}

\subsubsection{The relationship between the BIOs 
  and the solution operators}

\begin{theorem}\label{thm:1}
As an operator on $H^{1/2}(\Gamma)\times H^{-1/2}(\Gamma)$, $\First^{-1}$ has the decomposition
  \begin{equation}\label{eq:T1S1}
    \First^{-1} = S_{io} + S_{oi} -I
  \end{equation}
and,  as an operator on \emph{either} $H^{1/2}(\Gamma)\times H^{-1/2}(\Gamma)$ \emph{or} $H^1(\Gamma)\times L^2(\Gamma)$, $\Second^{-1}$ has the decomposition
  \begin{equation}\label{eq:T1S2}
    \Second^{-1} = I- S_{io} - S_{oi} + 2 S_{io}S_{oi}.
  \end{equation}
\end{theorem}

The proof of Theorem \ref{thm:1} is contained in \S\ref{sec:proof1} below.

The following result uses \eqref{eq:T1S1} and results about the behaviour of $S_{io}$ and
$S_{oi}$ in Lemmas \ref{lem:bound2} and \ref{lem:bound3} below to prove that if
$n_i\neq n_o$, then the norm of $\First$ blows up through the quasi-resonances of the transmission problem
\eqref{eq:BVP2} with $c_i = \max\{\nin,\nout\}$ and $c_o= \min\{\nin,\nout\}$.
This result explains rigorously the experiments in Figures \ref{fig:2d} and
\ref{fig:3d}. The result is stated using the weighted norm
$\N{\cdot}_{H^{1/2}_k(\Gamma)\times H^{-1/2}_k(\Gamma)}$ defined in \S\ref{sec:weighted},
with the operator norm
\begin{equation}\label{eq:abbreviation}
  \N{\cdot}_{H^{1/2}_k(\Gamma)\times H^{-1/2}_k(\Gamma)\rightarrow H^{1/2}_k(\Gamma)\times H^{-1/2}_k(\Gamma)}
  \quad\text{ abbreviated to } \quad
  \N{\cdot}_{H^{1/2}_k\times H^{-1/2}_k}.
\end{equation}

\begin{theorem}\mythmname{Superalgebraic blow up of $\|\First^{-1}\|$ for $\Oi$ smooth and convex}
  \label{thm:bound2}
  If $\Oi$ is $C^\infty$ with strictly-positive curvature and $n_i\neq n_o$, then there
  exist frequencies $0<k_1<k_2<\ldots$ with $k_j\rightarrow \infty$ such that given any
  $N>0$ there exists $C_N$ such that
  \begin{equation*}
    \N{\First^{-1}}_{H^{1/2}_{k_j}\times H^{-1/2}_{k_j}} \geq C_N k_j^N \quad\tfa j.
  \end{equation*}
\end{theorem}

The proof of Theorem \ref{thm:bound2} is contained in \S\ref{sec:proof2} below.

The reason we only prove blow up of $\First$, and not of $\Second$, is that Theorem \ref{thm:1} shows that $\Second^{-1}$ involves not only $S_{io}$ and $S_{oi}$ but also the composition of $S_{io}$ and $S_{oi}$ (whereas $\First$ does not), and we do not currently know how to
show that this extra term does not cancel out the blow up of one of  $S_{io}$ or $S_{oi}$.

The next result shows that, on appropriate subspaces, $\First^{-1}$ and $\Second^{-1}$ involve only the physical solution operator $S_{io}$. In particular, this result demonstrates that, because of the specific form of the right-hand sides in \eqref{eq:firstsecondkind}, only the physical solution operator $S_{io}$ is involved in the solution of the boundary value problem of Definition \ref{def:HTP}, as expected.
The results for $\Second^{-1}$ hold on either $H^{1/2}(\Gamma)\times H^{-1/2}(\Gamma)$ or 
$H^1(\Gamma)\times L^2(\Gamma)$, but the results for $\First^{-1}$ hold only on  $H^{1/2}(\Gamma)\times H^{-1/2}(\Gamma)$ (since we have not proved that $\First^{-1}$ exists on $H^1(\Gamma)\times L^2(\Gamma)$). We use the notation that $R(P)$ is the range of the operator $P$.

\begin{theorem}\mythmname{$\First$ and $\Second$ as operators $R(P_i^-) \rightarrow R(P_o^-)$}\label{thm:physical}

  (i) $\First^{-1} P_o^-=  \Second^{-1} P_o^-=S_{io} P_o^-$. 

  (ii) Both $\First$ and $\Second$ are bounded and invertible from $R(P_i^-)\rightarrow R(P_o^-)$ with $\First^{-1} = \Second^{-1}= S_{io}$ as operators from $R(P_o^-) \rightarrow R(P_i^-)$.
\end{theorem}

The proof of Theorem \ref{thm:physical} is contained in \S\ref{sec:proof1} below.

\subsubsection{Augmented BIEs}

We now propose a simple way to suppress spurious quasi-resonances in the BIEs without
  resorting to products of integral operators. We work in the Hilbert space
$\cH$ where $\cH:= H^{1/2}(\Gamma)\times H^{-1/2}(\Gamma)$ for the results involving
$\First$, and $\cH$ equals \emph{either} $H^{1/2}(\Gamma)\times H^{-1/2}(\Gamma)$ \emph{or}
$H^1(\Gamma)\times L^2(\Gamma)$ for the results involving $\Second$;
the norm $\|\cdot\|_{\cH}$ is then either
$\|\cdot\|_{H^{1/2}_k\times H^{-1/2}_k}$ or $\|\cdot\|_{H^{1}_k\times L^2}$.
We equip the space $\cH \times \cH$ with the norm
\begin{equation*}
  \N{\bpsi}_{\cH\times \cH}^2 := \N{\bpsi_1}_{\cH}^2 + \N{\bpsi_2}^2_{\cH},
\end{equation*}
where $\bpsi = (\bpsi_1, \bpsi_2)$ with $\bpsi_1, \bpsi_2 \in \cH$.

Define the augmented BIOs $\AugFirst$ and $\AugSecond : \cH \rightarrow \cH\times \cH$ by 
\begin{equation}\label{eq:Aug}
  \AugFirst := 
  \left(
    \begin{array}{c}
      \First \\
      P^+_i 
    \end{array}
  \right)
  \quad\tand
  \quad
  \AugSecond:=
  \left(
    \begin{array}{c}
      \Second \\
      P^+_i 
    \end{array}
  \right).
\end{equation}
The idea behind introducing these augmented operator equations is that the solution 
$\Cauchy^- u^-$ to the BIEs \eqref{eq:firstsecondkind} satisfies $P^+_i \Cauchy^- u^-=0$
(we see this below in \eqref{eq:L31} in the proof of Lemma \ref{lem:3}). 

\ble\mythmname{Solutions of augmented BIEs}
\label{lem:aug1}
Let $\widetilde{A}_*$ be one of $\AugFirst$ and $\AugSecond$.
Given ${\bf g} \in \cH$, if the solution $\bphi$ to the augmented operator equation
\begin{equation}\label{eq:augmented}
  \widetilde{A}_* \bphi = 
  \left(
    \begin{array}{c}
      {\bf g}\\
      {\bf 0}
    \end{array}
  \right)
\end{equation}
exists, then ${\bf g}$ satisfies
\begin{equation}\label{eq:constraint} 
  {\bf g}=S_{oi}{\bf g}
\end{equation}
and $\bphi$ is given by
\begin{equation}\label{eq:augsol}
  \bphi = S_{io}{\bf g}.
\end{equation}
\ele

Lemma \ref{lem:aug1} (proved in \S\ref{sec:proof_aug}) shows that the solution of the augmented operator equation \eqref{eq:augmented},
if it exists, only involves the physical solution operator $S_{io}$.  Note that if
${\bf g}= P_o^-{\bf f}$, i.e., the right-hand side of the first- and second-kind BIEs
\eqref{eq:firstsecondkind}, then \eqref{eq:constraint} is satisfied; indeed, it follows
from Lemma \ref{lem:SP} below that $(S_{oi}- I) P^-_o=0$.

\begin{example}
  \label{ex:aug}
  As in Example~\ref{ex:numexp} we perform a ``diagonalization'' of the augmented BIOs of
  \eqref{eq:Aug} to compute the operator norms of their pseudo-inverses in
  $H^{1/2}(\Gamma)\times H^{-1/2}(\Gamma)$ numerically (i.e., we compute the inverse
    of the smallest singular value of the block-diagonal matrix arising from truncating
    the Fourier/spherical-harmonic expansion). These norms as functions of the frequency
  $k$ are plotted in Figure \ref{fig:aug} for the case $n_{1}=1$, $n_{o}=3$, in which
  the physical solution operator $S_{io}$ has small norm for all values of $k$
  considered (as shown by the right-hand plots of Figures \ref{fig:2d} and \ref{fig:3d}).

  \begin{figure}[H]
\includegraphics[width=0.5\textwidth,clip]{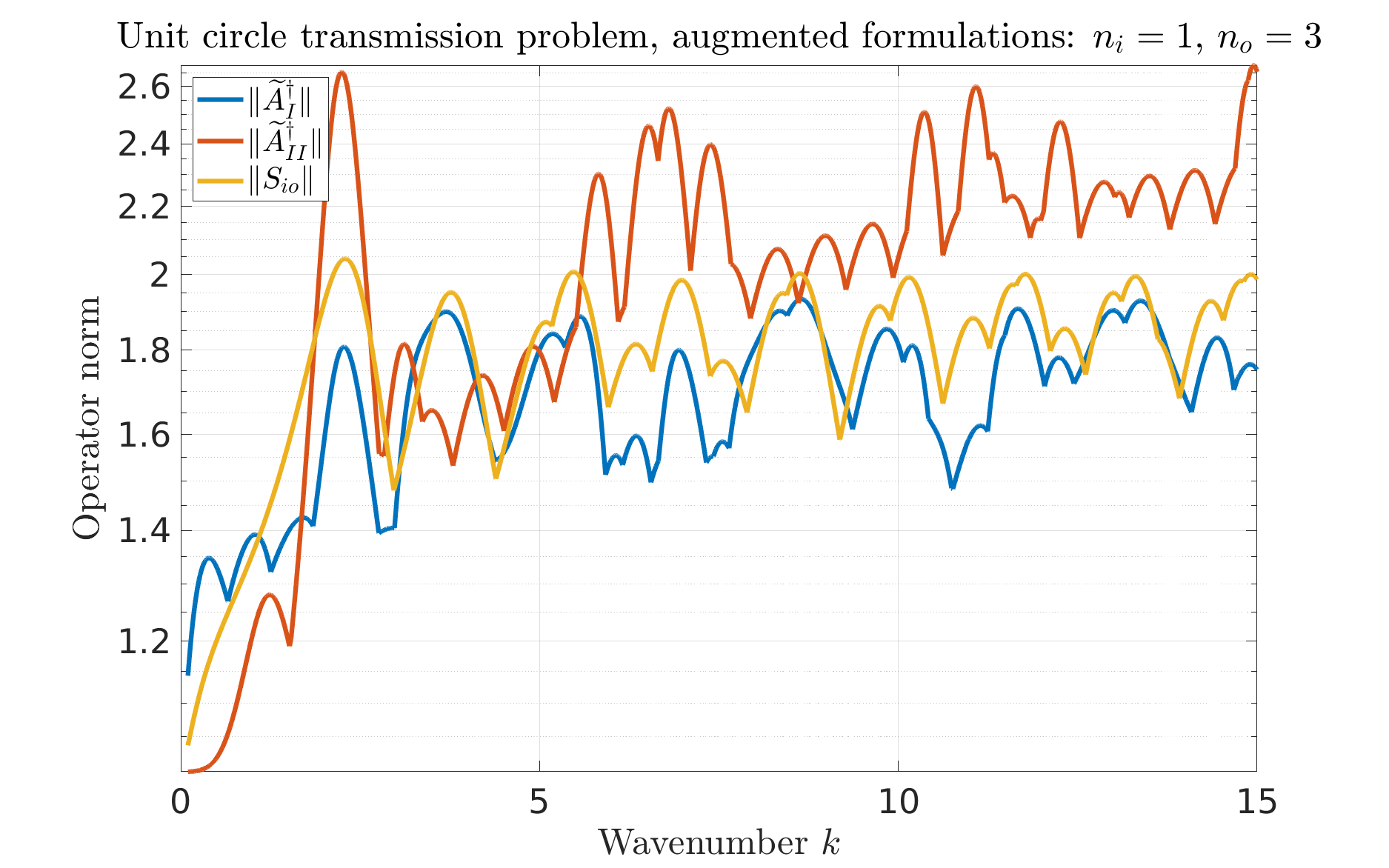}\hfill
\includegraphics[width=0.5\textwidth,clip]{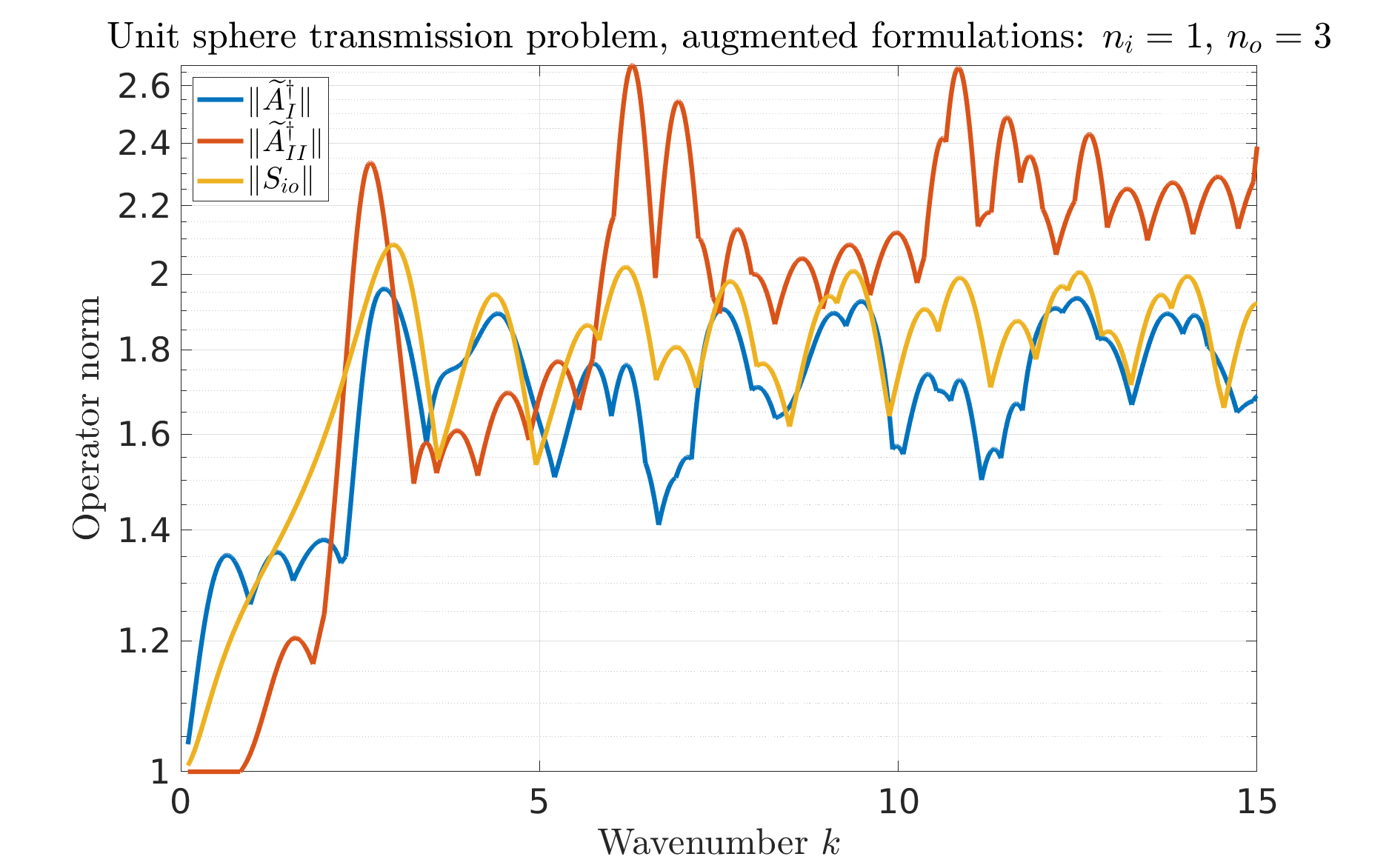}
\caption{Plots of the operator norms of the pseudo-inverses $\widetilde{A}_{\rm I}^{\dagger}$,
    $\widetilde{A}_{\rm II}^{\dagger}$ of the augmented BIOs}
    \label{fig:aug}
  \end{figure}

  As an agreeable surprise, we see that the norms of the pseudo-inverses of the augmented
  BIOs are smaller than those of $S_{io}$ for the range of frequencies considered -- augmentation has successfully removed any spurious quasi-resonances!
\end{example}

The following theorem rigorously explains the results in Figure \ref{fig:aug}, and is proved in \S\ref{sec:proof_aug}.

\begin{theorem}\mythmname{Stability of augmented BIEs}
  \label{thm:infsup}
  \begin{equation}\label{eq:infsup1}
    \inf_{\bphi\in\cH\setminus \{{\bf 0}\}}
    \sup_{\bpsi\in\cH\setminus \{{\bf 0}\}}
    \frac{
      \big|
      \big( \AugFirst \bphi, \bpsi\big)_{\cH\times\cH}
      \big|
    }{
      \N{\bphi}_{\cH} \N{\bpsi}_{\cH\times\cH}
    }
    \geq
    \frac{1}{
      \sqrt{2} \max\big\{ \N{S_{io}}_{\cH\rightarrow \cH}, 1\big\}
    }
  \end{equation}
  and
  \begin{equation}\label{eq:infsup1a}
    \inf_{\bphi\in\cH\setminus \{{\bf 0}\}}
    \sup_{\bpsi\in\cH\setminus \{{\bf 0}\}}
    \frac{
      \big|
      \big( \AugSecond \bphi, \bpsi\big)_{\cH\times\cH}
      \big|
    }{
      \N{\bphi}_{\cH} \N{\bpsi}_{\cH\times\cH}
    }
    \geq
    \frac{1}{
\sqrt{6+ 4\sqrt{2}}\, \max\big\{ \N{S_{io}}_{\cH\rightarrow \cH}, 1\big\}
    }.
  \end{equation}
\end{theorem}

This theorem reveals that the operator norms of the pseudo-inverses
$\widetilde{A}_{\rm I}^{\dagger}$ and $\widetilde{A}_{\rm II}^{\dagger}$ are bounded by
$C \max\big\{ \N{S_{io}}_{\cH\rightarrow \cH}, 1\big\}$ for some $k$-independent constant
$C>0$. Hence, if the physical solution operator $S_{io}$ is well-conditioned, then this well-conditioning carries over to the
BIOs of the augmented formulations.

\bre\textbf{(The analogue of Theorem \ref{thm:1} for BIOs for scattering by impenetrable
  obstacles.)}
\label{rem:impenetrable}
The analogous formulae to those in Theorem \ref{thm:1} for second-kind combined-field BIOs for solving the exterior
Dirichlet, Neumann, and impedance problems were given in \cite[Theorem 2.33]{CGLS12}, with
formulae for certain BIOs involving operator preconditioning given in \cite[Lemma
6.1]{BaSpWu:16}.  (We note that \cite[Lemma 6.1]{BaSpWu:16} introduced the idea of
obtaining these formulae via Calder\'on projectors, and we prove Theorem \ref{thm:1} using
this idea in \S\ref{sec:proof1}.)

For example, the standard direct second-kind combined-field BIO for solving the
exterior Dirichlet problem involves the operator $A_{\eta}':= \half I + K' - \ri \eta S$, for
$\eta\in \Rea\setminus\{0\}$, and \cite[Theorem 2.33]{CGLS12} and \cite[Lemma 6.1]{BaSpWu:16} (see also \cite[\S3]{EpGrHa:16}) prove that
\begin{equation}\label{eq:Euan_favourite_formula}
  \big(A'_{\eta}\big)^{-1}= I - \big({\rm DtN}^+ - \ri \eta \big){\rm ItD}^{-,\eta},
\end{equation}
where ${\rm DtN}^+$ is the exterior Dirichlet-to-Neumann map for solutions of the
Helmholtz equation satisfying the Sommerfeld radiation condition \eqref{eq:src}, and
${\rm ItD}^{-,\eta}$ is the interior Impedance-to-Dirichlet map (where the impedance
boundary condition is $\gamma_N^- u - \ri \eta \gamma^-_D u =g$). Recalling that $A'_\eta$
is also the standard indirect second-kind BIO for solving the interior impedance problem,
we see that \eqref{eq:Euan_favourite_formula} expresses $(A'_{\eta})^{-1}$ in terms of the
solution operators for the appropriate exterior and interior problems solved using
$A'_{\eta}$.  

The standard \emph{indirect} second-kind combined-field BIO for solving the
exterior Dirichlet problem involves the operator $A_{\eta}:= \half I + K - \ri \eta S$;
this operator is also the standard direct second-kind BIO for solving the interior impedance problem, and, correspondingly,
\begin{equation*}
  \big(A_{\eta}\big)^{-1}= I - {\rm ItD}^{-,\eta}\big({\rm DtN}^+ - \ri \eta \big).
\end{equation*}
\ere

\bre\textbf{(Indirect BIEs)} In this paper, we have considered only direct BIEs for the
Helmholtz transmission problem, i.e., BIEs where the unknown is the Cauchy data of the
solution. It is reasonable to expect that similar results hold for indirect BIEs for the 
transmission problem, just as similar decompositions into solution operators hold for the inverses of the direct BIOs for scattering by impenetrable
  obstacles (see the previous remark and \cite[Theorem 2.33]{CGLS12}), but we have not investigated this.
\ere

\bre\textbf{(Spurious quasi-resonances for electromagnetic BIEs)}
We expect that the phenomenon of spurious quasi-resonances also occurs for the BIEs for time-harmonic electromagnetic scattering; we have not pursued this in this paper however.
\ere

\section{Recap of results about layer potentials, BIOs, and Calder\'on projectors}\label{sec:background}

The single-layer and double-layer potentials, $\cV_{i/o}$ and $\cK_{i/o}$ respectively, are defined for $\varphi\in L^1(\Gamma)$ by 
\begin{align}\label{eq:SLP}
  \cV_{i/o} \varphi (\bx) &= \int_{\Gamma} \Phi_{i/o} (\bx,\by) \varphi (\by) \rd s (\by) \quad\tfa \bx \in \Rea^d \setminus \Gamma, \quad\tand \\
  \cK_{i/o} \varphi (\bx) &= \int_{\Gamma} \dfrac{\partial \Phi_{i/o} (\bx,\by)}{\partial n(\by)} \varphi (\by) \rd s (\by) \quad \tfa \bx \in \Rea^d \setminus \Gamma;
  \label{eq:DLP}
\end{align}
these definitions for $\varphi\in L^1(\Gamma)$ naturally extend to $\varphi \in H^{-s}(\Gamma)$ for $s\in[0,1]$ by continuity (see, e.g., \cite[Page 109]{CGLS12}).

\ble\label{lem:potentials}

(i) If $\phi\in H^{s-1/2}(\Gamma)$ with $|s|\leq 1/2$, then $\cV_{i/o} \phi\in H^{s+1}_{\rm loc}(\Rea^d)\cap C^2(\Rea^d\setminus\Gamma)\cap \SRC(k \sqrt{n_{i/o}})$.

(ii) If $\psi\in H^{s+1/2}(\Gamma)$ with $|s|\leq 1/2$, then $\cK_{i/o} \psi\in H^{s+1}_{\rm loc}(\Rea^d\setminus\Gamma)\cap C^2(\Rea^d\setminus\Gamma)\cap \SRC(k \sqrt{n_{i/o}})$.
\ele

\bpf[References for the proof]
See, e.g., \cite[Theorem 2.15]{CGLS12}; we note that the mapping properties for $|s|=1/2$ crucially use the harmonic analysis results of \cite{CoMcMe:82}, \cite{Ve:84}.
\epf

The potentials \eqref{eq:SLP} and \eqref{eq:DLP} are related to the integral operators in \eqref{eq:SD'} and \eqref{eq:DH} via the jump relations
\begin{equation}\label{eq:jumprelations}
  \gamma^\pm_D \cV_{i/o} = V_{i/o}, \quad \gamma^\pm_N \cV_{i/o} = \mp \frac{1}{2}I + K_{i/o}',
\quad  \gamma^\pm_D \cK_{i/o} = \pm \frac{1}{2}I + K_{i/o}, \quad\gamma^\pm_N \cK_{i/o} = -W_{i/o};
\end{equation}
see, e.g., \cite[\S7, Page 219]{MCL00}.
Recall the mapping properties, valid when $\Gamma$ is Lipschitz, $k\in \Com$, and $|s|\leq 1/2$,
 \begin{subequations}  \label{eq:mapping}
\begin{align}
  V_{i/o} : H^{s-1/2}(\Gamma)\rightarrow H^{s+1/2}(\Gamma), \quad\qquad& W_{i/o}: H^{s+1/2}(\Gamma)\rightarrow H^{s-1/2}(\Gamma),\\
  K_{i/o} : H^{s+1/2}(\Gamma)\rightarrow H^{s+1/2}(\Gamma), \quad\qquad& K'_{i/o}: H^{s-1/2}(\Gamma)\rightarrow H^{s-1/2}(\Gamma);
\end{align}
\end{subequations}
see, e.g., \cite[Theorems 2.17 and 2.18]{CGLS12} (similar to the results of Lemma \ref{lem:potentials}, the mapping properties for $|s|=1/2$ crucially use the harmonic analysis results of \cite{CoMcMe:82}, \cite{Ve:84}). 
The mapping properties \eqref{eq:mapping} imply that $P^{\pm}_{i/o}$ is a bounded operator from 
$H^{1/2}(\Gamma)\times H^{-1/2}(\Gamma)$ 
to itself and from 
$H^1(\Gamma)\times L^2(\Gamma)$ to itself.

We use the following notation for spaces of Helmholtz solutions: 
\begin{align*}
  \begin{aligned}
\Helmspace^-(\kappa):=&\big\{v\in H^1(\Om)\cap C^2(\Om),\; (\Delta+\kappa^2) v=0\big\},\\
\Helmspace^+(\kappa):=&\big\{v\in H^1\loc(\Op)\cap C^2(\Op)\cap\SRC(\kappa) ,\; (\Delta+\kappa^2) v=0\big\}.
  \end{aligned}
\end{align*}

\ble\label{lem:Cauchydata}
$R(P^{\pm}_{i/o})= \Cauchy^\pm \Helmspace^{\pm}(k\sqrt{n_{i/o}})$.
\ele

\bpf
By the jump relations \eqref{eq:jumprelations} and the definitions of $P^{\pm}_{i/o}$ \eqref{eq:Pdef}, with $\bphi= (\phi_1, \phi_2)$, 
\begin{equation}\label{eq:PBIE}
  P^{\pm}_{i/o} \bphi = \pm \Cauchy^{\pm}\big( \cK_{i/o}\phi_1 - \cV_{i/o}\phi_2\big);
\end{equation}
see, e.g., \cite[Equation 2.49]{CGLS12}. 
Both when $\bphi \in H^1(\Gamma)\times L^2(\Gamma)$ and when $\bphi \in H^{1/2}(\Gamma)\times H^{-1/2}(\Gamma)$, the right-hand side is then the trace of an element of $\Helmspace^{\pm}(k\sqrt{n_{i/o}})$ by Lemma \ref{lem:potentials}, so that $R(P^{\pm}_{i/o})\subset \Cauchy^\pm \Helmspace^{\pm}(k\sqrt{n_{i/o}})$.
To prove the reverse inclusion, given $u^\pm \in \Helmspace^\pm(k\sqrt{n_{i/o}})$, $u^\pm = \pm (\cK_{i/o} \gamma_D^\pm u - \cV_{i/o} \gamma_N^\pm u)$ by Green's integral representation (see, e.g., \cite[Theorems 2.20 and 2.21]{CGLS12});
\eqref{eq:PBIE} with $\phi_1 = \gamma^\pm_D u$ and $\phi_2 = \gamma_N^\pm u$ then implies that $\Helmspace^\pm(k\sqrt{n_{i/o}})\subset R(P^\pm_{i/o})$.
\epf

The following two lemmas are proved in, e.g., \cite[Page 118 and Lemma 2.22]{CGLS12}, respectively.\footnote{ Strictly speaking, \cite[\S2.5]{CGLS12} only considers $P^{\pm}_{i/o}$ as operators on $H^{1/2}(\Gamma)\times H^{-1/2}(\Gamma)$, but the proofs of the results on $H^1(\Gamma)\times L^2(\Gamma)$ are the same.}

\ble\label{lem:0}
$(P^+_{i/o})^2 = P^+_{i/o}$ and 
$(P^-_{i/o})^2 = P^-_{i/o}$ as operators \emph{either} on $H^{1/2}(\Gamma)\times H^{-1/2}(\Gamma)$ \emph{or} on $H^1(\Gamma)\times L^2(\Gamma)$.
\ele

\ble\label{lem:1}

(i) If $v\in\Helmspace^{-}(k\sqrt{n_{i/o}})$
then
\begin{equation}\label{eq:L11}
  P^-_{i/o} \Cauchy^- v = \Cauchy^- v.
\end{equation}

(ii) If $v\in\Helmspace^{+}(k\sqrt{n_{i/o}})$
then
\begin{equation}\label{eq:L12}
  P^+_{i/o} \Cauchy^+ v = \Cauchy^+ v.
\end{equation}
\ele

The next lemma is a converse to Lemma \ref{lem:1}.

\ble\label{lem:2}
Let $\bphi \in H^{1/2}(\Gamma)\times H^{-1/2}(\Gamma)$ \emph{or} $H^1(\Gamma)\times L^2(\Gamma)$.

(i) If $P^-_{i/o} \bphi =\bphi$, then $\bphi = \Cauchy^- v$ for
some $v\in\Helmspace^{-}(k\sqrt{n_{i/o}})$.

(ii) If $P^+_{i/o} \bphi =\bphi$, then $\bphi = \Cauchy^+ v$ for
some $v\in\Helmspace^{+}(k\sqrt{n_{i/o}})$.
\ele

\bpf
(i) 
Given $\bphi$ such that $P^-_{i/o}\bphi=\bphi$, let 
\begin{equation}\label{eq:v_-}
  v(\bx)= -\big( \cK_{i/o}\phi_1 - \cV_{i/o}\phi_2\big)(\bx) \quad\tfor \bx\in \Oi.
\end{equation}
By Lemma \ref{lem:potentials}, $v\in\Helmspace^{-}(k\sqrt{n_{i/o}})$.
We therefore only need to check that $\bphi = \Cauchy^- v$. However, by \eqref{eq:PBIE} and the definition of $v$ \eqref{eq:v_-},
$\bphi= P_{i/o}^- \bphi = \Cauchy^- v$.

(ii) 
Given $\bphi$ such that $P^+_{i/o}\bphi=\bphi$, let 
\begin{equation}\label{eq:v_+}
  v(\bx)= \big( \cK_{i/o}\phi_1 - \cV_{i/o}\phi_2\big)(\bx) \quad\tfor \bx\in \Oe.
\end{equation}
Similar to in (i), 
$v\in\Helmspace^{+}(k\sqrt{n_{i/o}})$,
and, by \eqref{eq:PBIE} and the definition of $v$ \eqref{eq:v_+}, 
$\bphi= P_{i/o}^+ \bphi = \Cauchy^+ v.$
\epf

\bpf[Proof of Lemma \ref{lem:3}]
By \eqref{eq:BVP} and \eqref{eq:L11}, $P_i^- \Cauchy^- u^- = \Cauchy^- u^-$, so that, by \eqref{eq:sum}, 
\begin{equation}\label{eq:L31}
  P^+_i \Cauchy^- u^-=0.
\end{equation}
Similarly, by \eqref{eq:BVP}, \eqref{eq:L12}, and \eqref{eq:sum},
\begin{equation}\label{eq:L32}
  P^-_o \Cauchy^+ u^+=0.
\end{equation}
Applying $P^-_o$ to the transmission condition $\Cauchy^-u^- = \Cauchy^+ u^+ + {\bf f}
$ in \eqref{eq:BVP} and using \eqref{eq:L32}, 
we  find that
\begin{equation}\label{eq:L34}
  P_o^- (\Cauchy^- u^-) = P_o^- {\bf f}.
\end{equation}

Subtracting \eqref{eq:L31} from \eqref{eq:L34}, we obtain the first-kind BIE in \eqref{eq:firstsecondkind}.
Adding \eqref{eq:L31} to \eqref{eq:L34}, we obtain the second-kind BIE in \eqref{eq:firstsecondkind}.
%

To obtain \eqref{eq:firstsecondkind2}, observe that, 
for the Helmholtz transmission scattering problem \eqref{eq:htsp}, ${\bf f} = \Cauchy^- u^I$ with $u^I$ satisfying \eqref{eq:entire}. Then $P_o^- \Cauchy^- u^I = \Cauchy^-u^I$ by \eqref{eq:L11}, and thus the right-hand sides of \eqref{eq:firstsecondkind} are just $\Cauchy^-u^I$.
\epf

\ble\label{lem:injective}
Let $A_*$ equal either $\First$ or $\Second$. Then $A_*$ is an injective, bounded operator on \emph{either} $H^{1/2}(\Gamma)\times H^{-1/2}(\Gamma)$ \emph{or} $H^1(\Gamma)\times L^2(\Gamma)$.

\ele

\bpf
The boundedness of $A_*$ follows from the expressions \eqref{eq:First}/\eqref{eq:Second} and the boundedness of $P_{i/o}^\pm$.
Injectivity follows by repeating the arguments in the proof of Theorem \ref{thm:1} below with ${\bf g=0}$ (these arguments use
uniqueness of the Helmholtz transmission problem of Definition \ref{def:HTP}).
\epf

\bpf[Proof of Lemma \ref{lem:invertible}]
The result for $\First$ follows from Lemma \ref{lem:injective} combined with the coercivity result in $H^{1/2}(\Gamma)\times H^{-1/2}(\Gamma)$ of, e.g., \cite[Theorem 7.27]{ClHiJePi:15} (see \cite[Corollary 7.28]{ClHiJePi:15}); we do not know of an analogous coercivity result in $H^1(\Gamma)\times L^2(\Gamma)$, hence why our results for $\First$ are only in $H^{1/2}(\Gamma)\times H^{-1/2}(\Gamma)$.

The results for $\Second$ follows from Lemma \ref{lem:injective} combined with the fact that $\Second- I$ is compact on both $H^{1/2}(\Gamma)\times H^{-1/2}(\Gamma)$ and $H^1(\Gamma)\times L^2(\Gamma)$. This latter result follows if we can show that 
\bit
\item $K_i-K_o$ is compact $H^{1/2}(\Gamma)\rightarrow H^{1/2}(\Gamma)$ and $H^{1}(\Gamma)\rightarrow H^1(\Gamma)$, 
\item $V_i-V_o$ is compact $H^{-1/2}(\Gamma)\rightarrow H^{1/2}(\Gamma)$ and $L^2(\Gamma)\rightarrow H^1(\Gamma)$,
\item $W_i-W_o$ is compact $H^{1/2}(\Gamma)\rightarrow H^{-1/2}(\Gamma)$ and $H^1(\Gamma)\rightarrow L^2(\Gamma)$, and
\item $K_i'-K_o'$ is compact $H^{-1/2}(\Gamma)\rightarrow H^{-1/2}(\Gamma)$ and $L^2(\Gamma)\rightarrow L^2(\Gamma)$.
  \eit
  Since $\Phi_{i}-\Phi_o =(\Phi_i-\Phi_0)-(\Phi_o-\Phi_0)$, where $\Phi_0$ is the Laplace fundamental solution, 
  these mapping properties follow from the bounds on 
  the difference of the Helmholtz and Laplace fundamental solutions in \cite[Equation 2.25]{CGLS12} and the fact that the inclusion $H^s(\Gamma)\rightarrow H^t(\Gamma)$ is compact for $-1\leq t\leq s\leq1$.
  \epf

  \section{Proof of Theorems \ref{thm:1} and \ref{thm:physical}}\label{sec:proof1}


  \ble\label{lem:4}
  Given $\bphi, {\bf f}$ in \emph{either} $H^{1/2}(\Gamma)\times H^{-1/2}(\Gamma)$ \emph{or} $H^{1}(\Gamma)\times L^2(\Gamma)$,
  \begin{equation}\label{eq:L41}
    \bphi = S_{io}{\bf f} \quad\text{ if and only if } 
    \quad
    \left\{
      \begin{array}{c}
        P_i^- \bphi=\bphi, \tand\\
        P_o^-(\bphi -{\bf f}) ={\bf 0}.
      \end{array}
    \right.
  \end{equation}
  Similarly,
  \begin{equation}\label{eq:L42}
    \bphi = S_{oi}{\bf f} \quad\text{ if and only if } 
    \quad
    \left\{
      \begin{array}{c}
        P_o^- \bphi=\bphi, \tand\\
        P_i^-(\bphi -{\bf f}) ={\bf 0}.
      \end{array}
    \right.
  \end{equation}
  \ele

  \bpf
  We prove \eqref{eq:L41}; the proof of \eqref{eq:L42} is the same with $i$ and $o$ swapped.

  We first prove the forward implication in \eqref{eq:L41}. Given ${\bf f}$, let $u$ be as in the definition of $S_{io}$ (Definition \ref{def:Sio}), i.e., $u$ satisfies \eqref{eq:BVP2} with $c_i=n_i$ and $c_o=n_o$. By definition $\bphi =\Cauchy^- u$, so $P_i^-\bphi=\bphi$ by \eqref{eq:L11}. The jump condition in \eqref{eq:BVP2} implies that $\bphi- {\bf f}=\Cauchy^+ u$, and \eqref{eq:L12} then implies that $P^+_o(\bphi-{\bf f}) = \bphi-{\bf f}$.

  For the reverse implication in \eqref{eq:L41}, given $\bphi$ satisfying the right-hand side of \eqref{eq:L41}, 
  Part (i) of Lemma \ref{lem:2} implies that $\bphi = \Cauchy^- w^-$ for some $w^-\in \Helmspace^{-}(k\sqrt\nin)$.
  Similarly, Part (ii) of Lemma \ref{lem:2} implies that 
  $\bphi -{\bf f}= \Cauchy^+ w^+$ for some
  $w^+\in \Helmspace^{+}(k\sqrt\nout)$.
  Let  $w:= w^+$ in $\Oe$ and $w:= w^-$ in $\Oi$.
  Then $\Cauchy^- w^- - \Cauchy^+ w^+ = \bphi- (\bphi-{\bf f})= {\bf f}$. Since the solution of the transmission problem is unique, $w$ equals the function $u$ in the definition of $S_{io}$ (i.e., Definition \ref{def:Sio} with $c_i=n_i$ and $c_o=n_o$), and $\bphi = \Cauchy^-w^- =\Cauchy^-u^- = S_{io}{\bf f}$.
  \epf

  We now prove Theorem \ref{thm:1}.

  \bpf[Proof of the result \eqref{eq:T1S1} in Theorem \ref{thm:1}]
  Assume that $\bpsi, {\bf g} \in H^{1/2}(\Gamma)\times H^{-1/2}(\Gamma)$ or $H^{1}(\Gamma)\times L^2(\Gamma)$ with 
  $\First \bpsi= {\bf g}$, i.e.,
  \begin{equation}\label{eq:T11}
    (P^-_i - P^+_o)\bpsi = {\bf g}.
  \end{equation}

  \noindent\emph{Step 1: apply $P^-_i$ to \eqref{eq:T11}}. 

  \noindent Applying $P^-_i$ to \eqref{eq:T11} and using the fact that  $P^-_i$ is a projection (by Lemma \ref{lem:0}), we have
  \begin{equation*}
    P_i^- \big( \bpsi - P_o^+ \bpsi - {\bf g}\big) = {\bf 0},
  \end{equation*}
  that is, by \eqref{eq:sum},
  \begin{equation}\label{eq:T11a}
    P_i^- \big( P_o^-\bpsi- {\bf g}\big) = {\bf 0},
  \end{equation}
  Let $\bphi := P_o^-\bpsi$ and let ${\bf f :=g}$. Then by Lemma \ref{lem:4} and the fact that  $P^-_o$ is a projection, $\bphi= S_{oi} {\bf g}$, i.e. 
  \begin{equation}\label{eq:T12}
    P_o^- \bpsi = S_{oi}{\bf g}.
  \end{equation}

  \noindent\emph{Step 2: apply $P^-_o$ to \eqref{eq:T11}}. 

  \noindent Applying $P^-_o$ to \eqref{eq:T11} and using the fact that  $P^-_o$ is a projection (so that, in particular, $P^-_oP^+_o=0$), we have
  \begin{equation}\label{eq:T12a}
    P_o^- \big( P_i^- \bpsi - {\bf g}\big) = {\bf 0}.
  \end{equation}
  Let $\bphi := P_i^-\bpsi$ and let ${\bf f :=g}$. Then by Lemma \ref{lem:4} and the fact that  $P^-_i$ is a projection, $\bphi= S_{io} {\bf g}$, i.e.
  \begin{equation}\label{eq:T13}
    P_i^- \bpsi = S_{io}{\bf g}.
  \end{equation}

  \noindent \emph{Step 3: use \eqref{eq:sum} and \eqref{eq:T11} and the results of Steps 1 and 2}. 

  \noindent By \eqref{eq:sum},  \eqref{eq:T12}, \eqref{eq:T11}, and \eqref{eq:T13} (in that order),
  \begin{align}\label{eq:end1}
    \bpsi = (P_o^- + P_o^+) \bpsi &= S_{oi}{\bf g} + P_i^- \bpsi - {\bf g} = (S_{oi} + S_{io} -I){\bf g},
  \end{align} 
  which is the result \eqref{eq:T1S1}.
  \epf

  \bpf[Proof of the result \eqref{eq:T1S2} in Theorem \ref{thm:1}]
  Assume that $\bpsi, {\bf g} \in H^{1/2}(\Gamma)\times H^{-1/2}(\Gamma)$ 
  or $H^{1}(\Gamma)\times L^2(\Gamma)$
  with 
  $\Second \bpsi= {\bf g}$, i.e.,
  \begin{equation}\label{eq:T14}
    (P^-_o + P^+_i)\bpsi = {\bf g}.
  \end{equation}

  \noindent\emph{Step 1: apply $P^-_i$ to \eqref{eq:T14}}. 

  \noindent Applying $P^-_i$ to \eqref{eq:T14} and using the fact that  $P^-_i$ is a projection (by Lemma \ref{lem:0}), we see that \eqref{eq:T11a} holds.
  Let $\bphi := P_o^-\bpsi$ and let ${\bf f :=g}$. Then by Lemma \ref{lem:4} and the fact that  $P^-_o$ is a projection, $\bphi= S_{oi} {\bf g}$, i.e., \eqref{eq:T12} holds.

  \

  \noindent\emph{Step 2: apply $P^-_o$ to \eqref{eq:T14}}. 

  \noindent Applying $P^-_o$ to \eqref{eq:T14} and using the fact that  $P^-_o$ is a projection, we see that
  \begin{equation*}
    P_o^- \big(\bpsi + P_i^+ \bpsi - {\bf g}\big) = {\bf 0}.
  \end{equation*}
  Let $\bphit := P_i^-\bpsi$, so that 
  \begin{equation*}
    P_o^-\big(\bphit +2 P_i^+\bpsi -{\bf g}\big) ={\bf 0}.
  \end{equation*}
  Let ${\bf \widetilde{f}} :=-2 P_i^+\bpsi +{\bf g}$. Then by Lemma \ref{lem:4} and the fact that  $P^-_i$ is a projection, 
  $\bphit = S_{io} {\bf \widetilde f}$, i.e.,
  \begin{equation*}
    P_i^- \bpsi= S_{io}\big(-2 P_i^+\bpsi +{\bf g}\big).
  \end{equation*}
  Using \eqref{eq:T14} and then \eqref{eq:T12}, which holds by \emph{Step 1}, we have
  \begin{equation}\label{eq:T16}
    P_i^- \bpsi= S_{io}\big(2 P_o^-\bpsi -{\bf g}\big)
    =S_{io}\big(2 S_{oi} {\bf g}-{\bf g}\big).
  \end{equation}

  \noindent \emph{Step 3: use \eqref{eq:sum}, \eqref{eq:T14}, and the results of Steps 1 and 2}. 

  \noindent By \eqref{eq:sum} and \eqref{eq:T12}, 
  \begin{align}\label{eq:T17}
    \bpsi = (P_o^- + P_o^+) \bpsi = S_{oi}{\bf g} + P_o^+ \bpsi.
  \end{align} 
  Using \eqref{eq:sum} in \eqref{eq:T14} and rearranging, we have
  \begin{equation*}
    P_o^+\bpsi = - P_i^- \bpsi + 2\bpsi -{\bf g},
  \end{equation*}
  and using this in \eqref{eq:T17} we find that 
  \begin{equation*}
    \bpsi  =\big(I- S_{oi}\big){\bf g} + P_i^-\bpsi;
  \end{equation*}
  the result \eqref{eq:T1S2} then follows from using \eqref{eq:T16}.
  \epf

  To prove Theorem \ref{thm:physical}, we need the following consequences of the definitions of $S_{io}$ and $S_{oi}$.

  \ble\label{lem:SP}
  $S_{io}P^+_o=0$ and $S_{io}P^-_i=P_i^-$
  as operators on \emph{either} $H^{1/2}(\Gamma)\times H^{-1/2}(\Gamma)$ \emph{or} $H^1(\Gamma)\times L^2(\Gamma)$.
  Similarly, $S_{oi}P^+_i=0$ and $S_{oi}P^-_o=P^-_o$.
  \ele

  \bpf
  We prove the relationships involving $S_{io}$; the proofs of those involving $S_{oi}$ are completely analogous.
  Given ${\bf f }\in H^{1/2}(\Gamma)\times H^{-1/2}(\Gamma)$ or $H^1(\Gamma)\times L^2(\Gamma)$, by Definition \ref{def:Sio}, $S_{io} P^+_o {\bf f} = \Cauchy^- v^-$ where 
  \begin{equation}\label{eq:Thursday1}
    v^+ \in \Helmspace^{+}(k\sqrt{n_o}), \quad v^- \in \Helmspace^{-}(k\sqrt{n_i}), \quad \tand \quad \Cauchy^- v^- = \Cauchy^+ v^+ + P^+_o{\bf f}.
  \end{equation}
  By Lemma \ref{lem:Cauchydata} and \eqref{eq:PBIE}, there exists $w^+ \in \Helmspace^{+}(k\sqrt{n_o})$ such that $\Cauchy^+ w^+ = P_o^+{\bf f}$. Thus $v^-:=0$ and $v^+ := - w^+$ is a solution of \eqref{eq:Thursday1}, and by uniqueness of the Helmholtz transmission problem (Lemma \ref{lem:wellposed}) it is the only solution. Therefore $S_{io} P^+_o {\bf f} =\Cauchy^- v^-={\bf 0}$. 

  The proof that $S_{io} P^-_i=P^-_i$ is similar. Indeed, again using uniqueness of the Helmholtz transmission problem, we have $S_{io} P^-_i {\bf f} = \Cauchy^- w^-$ with $w^- \in \Helmspace^{-}(k\sqrt{n_i})$ and $\Cauchy^- w^- = P_i^- {\bf f}$.
  \epf

  %

  \bpf[Proof of Theorem \ref{thm:physical}]
  Part (i): By \eqref{eq:T1S1},
  $\First^{-1}P_o^- = (S_{io}+ S_{oi}-I) P_o^-.$
  Lemma \ref{lem:SP} shows that $(S_{oi}-I)P^-_o=0$, and thus $\First^{-1} P_o^-=S_{io} P_o^-$. The equation $\Second^{-1} P_o^-=S_{io} P_o^-$ follows similarly.

  Part (ii): By the second equality in \eqref{eq:First}, Lemma \ref{lem:0}, and \eqref{eq:sum},
  \begin{equation*}
    \First P_i^- = (P_i^-- P_o^+)P_i^-= (I - P_o^+) P_i^- = P_o^- P_i^-,
  \end{equation*}
  so that $\First: R(P_i^-) \rightarrow R(P_o^-)$. 
  Similarly, by \eqref{eq:Second},
  \begin{equation*}
    \Second P_i^- = (P_o^- + P_i^+) P_i^- = P_o^- P_i^-,
  \end{equation*}
  so that $\Second: R(P_i^-) \rightarrow R(P_o^-)$. 
  By Definition \ref{def:Sio}, $S_{io}$ maps into the space of Cauchy data of $\Helmspace^{-}(k\sqrt{n_i})$; by 
  \eqref{eq:L11} and Lemma \ref{lem:Cauchydata}, this space is $R(P_i^-)$. Therefore, by Part (i), both $\First^{-1}$ and $\Second^{-1}$ map
  $R(P_o^-) \rightarrow R(P_i^-)$ and both equal $S_{io}$ as operators between these spaces.
  \epf

  \section{Proof of Theorem \ref{thm:bound2}}\label{sec:proof2}

  Throughout this section we use the notation that $a\lesssim b$ if there exists $C>0$, independent of $k$, such that $a\leq Cb$. We write 
  $a\sim b$ if both $a\lesssim b$ and $b\lesssim a$.

  \subsection{Definitions of \texorpdfstring{$k$}{k}-weighted norms and associated results}\label{sec:weighted}
  For $\bphi \in H^1(\Gamma)\times L^2(\Gamma)$ with $\bphi= (\phi_1,\phi_2)$, let 
  $\nabla_T\phi_1$ be the tangential gradient of $\phi_1$ on $\Gamma$ and
  \begin{equation*}
    \N{\phi_1}_{H^1_k(\Gamma)}^2:= \N{\nabla_T \phi_1}_{L^2(\Gamma)}^2 +k^2\N{\phi_1}_{L^2(\Gamma)}^2,
    \qquad
    \N{\bphi}_{H^1_k(\Gamma)\times L^2(\Gamma)}^{2}:= 
    \N{\phi_1}_{H^1_k(\Gamma)}^{2} + \N{\phi_2}_{L^2(\Gamma)}^{2}.
  \end{equation*}
  Define $H^{1/2}_k(\Gamma)$ by interpolation between $H^1_k(\Gamma)$ and $L^2(\Gamma)$ and then $H^{-1/2}_k(\Gamma)$ by duality. As in \S\ref{sec:intro}, we use the abbreviation \eqref{eq:abbreviation}.

  For a bounded Lipschitz open set $D\subset \Rea^d$, let
  \begin{equation*}
    \N{v}_{H^1_k(D)}^2:= \N{\nabla v}^2_{L^2(D)}+k^2 \N{v}^2_{L^2(D)}.
  \end{equation*}
Fix $k_0>0$. Then, with $H^{1/2}_k(\partial D)$ defined above, by, e.g., \cite[Theorem 5.6.4]{NED01},
  \begin{equation}\label{eq:Dtrace}
    \N{\gamma_D v}_{H^{1/2}_k(\partial D)}\lesssim \N{v}_{H^1_k(D)} \quad\tfa v\in H^1(D) \;\tand k\geq k_0,
  \end{equation}
  and there exists $E: H^{1/2}(\partial D)\to H^1(D)$ such that
  \begin{equation}\label{eq:extension}
    \gamma_DE\phi=\phi\quad\tand\quad
    \N{E\phi}_{H^1_k(D)}\lesssim \N{\phi}_{H^{1/2}_k(\partial D)}.
  \end{equation}


  \ble\label{lem:Neumann}
  If $v\in H^1(D,\Delta)$ with $(\Delta + k^2 c)v=0$ and $k\geq k_0$, then
  \begin{equation}\label{eq:Ntrace}
    \N{\gamma_N v}_{H^{-1/2}_k(\partial D)}\lesssim \N{v}_{H^1_k(D)}
  \end{equation}
  (where the omitted constant depends on $c$ and $k_0$).
  \ele

  \bpf[Sketch proof of Lemma \ref{lem:Neumann}]
  This follows by repeating the argument in, e.g., \cite[Lemma 4.3]{MCL00} (which starts from the definition of the Neumann trace via Green's identity) and then using weighted norms and, in particular, the bound \eqref{eq:extension}.
  \epf

  \subsection{From resolvent estimates to bounds on \texorpdfstring{$S_{io}$ and $S_{oi}$}{Sio and Soi}}


  \ble\label{lem:bound1}
  Given $c_o, c_i, R$ positive real numbers and ${\bf g} \in H^{1/2}(\Gamma)\times H^{-1/2}(\Gamma)$, let $v \in H^1_{\rm loc}(\Rea^d\setminus \Gamma)\cap \SRC( k\sqrt{c_o})$ satisfy
  \begin{align*}
    \begin{aligned}
      (\Delta +k^2 c_i)v^- &=0&&\iin\Oin,\\
      (\Delta +k^2 c_o)v^+ &=0 &&\iin\Oout,\\
      \Cauchy^- v^- &= \Cauchy^+ v^++ {\bf g} &&\oon \Gamma.
    \end{aligned}
  \end{align*}
  Assume that, for all $f\in L^2(\Rea^d)$ with $\supp f \subset B_R$ and $w \in H^1_{\rm loc}(\Rea^d\setminus \Gamma)\cap \SRC( k\sqrt{c_o})$ that satisfy
  \begin{align*}
    \begin{aligned}
      (\Delta +k^2 c_i)w^- &=f^-&&\iin\Oin,\\
      (\Delta +k^2 c_o)w^+ &=f^+ &&\iin\Oout,\\
      \Cauchy^- w^- &= \Cauchy^+ w^+  &&\oon \Gamma,
    \end{aligned}
  \end{align*}
  the following bound holds:
  \begin{equation}\label{eq:Csol}
    \N{w}_{H^1_k(B_R)}\leq C_{\rm sol}(k, R, c_i,c_o) \N{f}_{L^2(B_R)}.
  \end{equation}
  Then, given $k_0>0$,
  \begin{equation}\label{eq:Csol2}
    \N{v}_{H^1_k(B_R)}\lesssim k\big(1 + C_{\rm sol}(k, R, c_i,c_o)\big) \N{g}_{H^{1/2}_k(\Gamma)\times H^{-1/2}_k(\Gamma)} \quad\tfa k \geq k_0.
  \end{equation}
  \ele

  \begin{cor}\label{cor:bound1} 
    Under the assumptions of Lemma \ref{lem:bound1},
    \begin{equation*}
      \N{S(c_i,c_o)}_{H^{1/2}_k\times H^{-1/2}_k}
      \lesssim k\big(1 + C_{\rm sol}(k, R, c_i,c_o)\big).
    \end{equation*}
  \end{cor}

  \bpf[Proof of Corollary \ref{cor:bound1} from Lemma \ref{lem:bound1}]
  This follows from combining the result of Lemma \ref{lem:bound1} and the trace results \eqref{eq:Dtrace} and \eqref{eq:Ntrace}.
  \epf

  \bpf[Proof of Lemma \ref{lem:bound1}]
  Let $u\in H^1(\Rea^d\setminus \Gamma)\cap \SRC(k\sqrt{c_o})$ be the solution to 
  \begin{align}
    \begin{aligned}
      (\Delta + (k^2+\ri k)c_i)u^- &=0&&\iin\Oin,\\
      (\Delta + (k^2+\ri k)c_o)u^+ &=0 &&\iin\Oout\cap B_{R'},\\
      \Cauchy^- u^- &= \Cauchy^+ u^+ + {\bf g} &&\oon \Gamma
    \end{aligned}
    \label{eq:absorption}
  \end{align}
  (this choice of auxiliary problem is motivated by the proof of \cite[Theorem 3.5]{BaSpWu:16} using \cite[Lemma 3.3]{BaSpWu:16}). We prove below that 
  \begin{equation}\label{eq:absorption_bound}
    \N{u}_{H^1_k(\Oout)} \lesssim k\N{{\bf g}}_{H^{1/2}_k(\Gamma)\times H^{-1/2}_k(\Gamma)}.
  \end{equation}
Given $R>0$ such that $\Oi \Subset B_R$, choose $\chi \in C^\infty_{\rm comp}(\Rea^d)$ with $\supp \chi \subset B_R$ and $\chi\equiv1$ on $\Oin$. 
  Let $w=v- \chi u$; then $w\in 
  H^1_{\rm loc}(\Rea^d\setminus \Gamma)\cap \SRC(k\sqrt{c_o})$ satisfies 
  \begin{align*}
    \begin{aligned}
      (\Delta +k^2 c_i)w^- &=\ri k c_i u^-&&\iin\Oin,\\
      (\Delta +k^2 c_o)w^+ &=\ri k c_o \chi u^+ - 2 \nabla u^+ \cdot\nabla \chi - u^+\Delta \chi &&\iin\Oout,\\
      \Cauchy^- w^- &= \Cauchy^+ w^+  &&\oon \Gamma.\\
    \end{aligned}
  \end{align*}
  Using the fact that $w=v-\chi u$, the fact that $\supp \chi \subset B_R$ (by construction), and the bound \eqref{eq:Csol}, we have that, given $k_0>0$,
  \begin{align*}
    \N{v}_{H^1_k(B_R)} &\lesssim \N{w}_{H^1_k(B_R)} + \N{u}_{H^1_k(B_R)} \lesssim C_{\rm sol}(k, R, c_i,c_o)  \N{u}_{H^1_k(B_R)} +  \N{u}_{H^1_k(B_R)}
  \end{align*}
for all $k\geq k_0$; the result \eqref{eq:Csol2} then follows from the bound \eqref{eq:absorption_bound}.

  It therefore remains to prove \eqref{eq:absorption_bound}. First observe that, thanks to the $\ri k$ term in the PDE in \eqref{eq:absorption}, $u^+$ decays exponentially at infinity, and thus $u^+ \in H^1(\Oout)$. Next, apply Green's identity to $u^-$ in $\Oin$ and $u^+$ in $\Oout $ to obtain that
  \begin{align}\label{eq:Green1}
    \N{\nabla u^-}^2_{L^2(\Oin)} - (k^2+ \ri k) c_i \N{u^-}^2_{L^2(\Oin)} &= \langle \gamma_N^- u^-, \gamma_D^- u^-\rangle_\Gamma, \\
    \N{\nabla u^+}^2_{L^2(\Oout)} - (k^2+ \ri k) c_o \N{u^+}^2_{L^2(\Oout)} &= -\langle \gamma_N^+ u^+, \gamma_D^+ u^+\rangle_\Gamma.
    \label{eq:Green2}
  \end{align}
  The jump condition in \eqref{eq:absorption} implies that, with ${\bf g}= (g_D, g_N)$, 
  \begin{equation}\label{eq:Green2a}
    \langle \gamma_N^- u^-, \gamma_D^- u^-\rangle_\Gamma- \langle \gamma_N^+ u^+, \gamma_D^+ u^+\rangle_\Gamma
    = \langle g_N, \gamma_D^- u^- \rangle_\Gamma + \langle \gamma_N^+ u , g_D\rangle_\Gamma.
  \end{equation}
  Therefore, adding \eqref{eq:Green1} and \eqref{eq:Green2}, taking the imaginary part, and then using the Cauchy-Schwarz inequality on terms arising from \eqref{eq:Green2a}, we obtain
  \begin{equation}\label{eq:Green3}
    \min\{c_i,c_o\} k\N{u}^2_{L^2(\Oout)} \leq \N{g_N}_{H^{-1/2}_k(\Gamma)}\N{\gamma_D^- u^-}_{H^{1/2}_k(\Gamma)}
    +\N{\gamma_N^+ u^+}_{H^{-1/2}_k(\Gamma)}\N{g_D}_{H^{1/2}_k(\Gamma)}.
  \end{equation}
  Adding \eqref{eq:Green1} and \eqref{eq:Green2}, taking the real part, adding a sufficiently-large multiple of $k$ times \eqref{eq:Green3}, 
  and then using the Cauchy-Schwarz inequality on terms arising from \eqref{eq:Green2a}, 
  we have 
  \begin{equation}\label{eq:Green4}
    \N{u}^2_{H^1_k(\Oout)} \lesssim k \Big( \N{g_N}_{H^{-1/2}_k(\Gamma)}\N{\gamma_D^- u^-}_{H^{1/2}_k(\Gamma)}
    +\N{\gamma_N^+ u^+}_{H^{-1/2}_k(\Gamma)}\N{g_D}_{H^{1/2}_k(\Gamma)}\Big).
  \end{equation}
  The bound \eqref{eq:absorption_bound} then follows from using the inequality
  \begin{equation}\label{eq:Cauchy}
    2ab \leq \epsilon a^2 + \epsilon^{-1} b^2, \quad a,b,\epsilon>0,
  \end{equation}
  and the trace bounds \eqref{eq:Dtrace} and \eqref{eq:Ntrace} in the right-hand side of \eqref{eq:Green4}.
  \epf

  \subsection{\texorpdfstring{$k$-explicit bounds on $S_{io}$ and $S_{oi}$}{k-explicit bounds on Sio and Soi}}


  We recall the notions of  \emph{star-shaped} and \emph{star-shaped with respect to a ball}. 

  \begin{defin}
    (i) $\Oin$ is \emph{star-shaped with respect to the point $\bx_0$} if, whenever $\bx \in \Oin$, the segment $[\bx_0,\bx]\subset \Oin$.

    \noindent (ii) $\Oin$ is \emph{star-shaped with respect to the ball $B_{a}(\bx_0)$} if it is star-shaped with respect to every point in $B_{a}(\bx_0)$.
  \end{defin}

  \begin{lemma}\mythmname{``Good'' behaviour of $S_{io}$ when $n_i<n_o$}\label{lem:bound2}
    If $\Oi$ is star-shaped with respect to a ball and $n_i<n_o$, then, given $k_0>0$, 
    \begin{equation*}
      \N{S_{io}}_{H^{1/2}_k\times H^{-1/2}_k}
       \lesssim k\quad\tfa k\geq k_0.
    \end{equation*}
  \end{lemma}

  \bpf
  \cite[Theorem 3.2]{MoSp:19} proves that \eqref{eq:Csol} holds with $C_{\rm sol}(k)\sim 1$, and the result then follows from Corollary \ref{cor:bound1}.
  \epf

  %
  %
  %

  \ble\mythmname{``Bad'' behaviour of $S_{io}$ when $n_i>n_o$}\label{lem:bound3}
  If $\Oi$ is $C^\infty$ with strictly-positive curvature and $n_i>n_o$, then there exist $0<k_1<k_2<\ldots$ with $k_j\rightarrow \infty$ such that given any $N>0$ there exists $C_N>0$ such that 
  \begin{equation*}
    \N{S_{io}}_{H^{1/2}_k\times H^{-1/2}_k} \geq C_N k_j^N \quad\tfa j.
  \end{equation*}
  \ele

  In the proof of Lemma \ref{lem:bound3}, we use the notation that $a = O(k^{-\infty})$ as $k\tendi$ if, given $N>0$, there exists $C_N,k_0$ such that $|a|\leq C_N k^{-N}$ for all $k\geq k_0$, i.e.~$a$ decreases superalgebraically in $k$. 

  The ideas behind Lemma \ref{lem:bound3} are that (i) if there exist  quasimodes with $O(k^{-\infty})$ remainder (in the sense of  \eqref{eq:QM1} below), then the norm of $S_{io}$ has $O(k^{\infty})$ blow up (immediately from the definitions of quasimodes and $S_{io}$), and (ii)  if $\Oi$ is $C^\infty$ with strictly positive curvature and $n_i>n_o$ then  quasimodes with $O(k^{-\infty})$ remainder exist by \cite{PoVo:99}.
  To prove Lemma \ref{lem:bound3}, we need the following bounds on the Newtonian potential, i.e., integration against the fundamental solution. 
  Let
  \begin{align*}
    \cN_{i/o} f (x):= \int_{\Rea^d} \Phi_{i/o} (x,y)f(y) \, \rd y.
  \end{align*}

  \begin{lemma}
    \label{lem:Newton}
    Given $f\in L^2(\Rea^d)$ with $\supp f \subset B_R$ and $k_0>0$,
    \begin{equation*}
      \frac{1}{k}
  \sum_{|\alpha|=2}\N{\partial^\alpha (\cN_{i/o} f)}_{L^2(B_R)}
      + 
      \N{\cN_{i/o} f}_{H^1_k(B_R)}
      \lesssim \N{f}_{L^2(B_R)}
    \end{equation*}
    for all $k\geq k_0$, where the omitted constant depends on $n_{i/o}$ and $R$.
  \end{lemma}

  \bpf[References for the proof of Lemma \ref{lem:Newton}]
  See, e.g., \cite[Theorem 3.1]{DyZw:19} for $d=3$ and \cite[Theorem 14.3.7]{Ho:83a} for arbitrary dimension (note that \cite[Theorem 14.3.7]{Ho:83a} is for fixed $k$, but a rescaling of the independent variable yields the result for arbitrary $k$).
  \epf

  \bpf[Proof of Lemma \ref{lem:bound3}]
  We are going to show that there exist ${\bf g}_j \in H^{1/2}(\Gamma)\times H^{-1/2}(\Gamma)$, $j=1,2,\ldots$, such that the solutions $v_j$ to \eqref{eq:BVP2} (with $c_i=n_i$ and $c_o=n_o$ and ${\bf f= g}_j$) are such that given $N>0$ there exists $C_N>0$ such that 
  \begin{equation*}%
    \N{\Cauchy^- v_j^-}_{H^{1/2}_k(\Gamma)\times H^{-1/2}_k(\Gamma)} \geq C_N k_j^N 
    \N{{\bf g}_j}_{H^{1/2}_k(\Gamma)\times H^{-1/2}_k(\Gamma)}
    \quad\tfa j;
  \end{equation*}
  i.e., that 
  \begin{equation}\label{eq:QM0}
    \frac{
      \N{{\bf g}_j}_{H^{1/2}_k(\Gamma)\times H^{-1/2}_k(\Gamma)}
    }{
      \N{\Cauchy^- v_j^-}_{H^{1/2}_k(\Gamma)\times H^{-1/2}_k(\Gamma)}
    } = O(k_j^{-\infty})\quad\tas j\tendi.
  \end{equation}

  By \cite{PoVo:99}, there exist $k_j \in \Com$, with $|k_j|\tendi$, $0>\Im k_j = O(k_j^{-\infty})$, 
  $w_j^{\pm}\in C^\infty(\overline{\Omega^{\pm}})$ with the support of $w_j^{\pm}$ contained in a fixed compact neighbourhood of $\Gamma$, and such that $\|\gamma_D^- w^-_j\|_{L^2(\Gamma)}=1$, 
  \begin{align}\label{eq:QM1}
    \N{(\Delta + k_j^2 n_{i/o}) w_j^{\pm}}_{L^2(\Omega^{\pm})} = O(k_j^{-\infty}), \qquad
    \N{\Cauchy^- w_j^- - \Cauchy^+ w_j^+}_{H^{2}(\Gamma)\times H^{2}(\Gamma)} = O(k_j^{-\infty}),
  \end{align}
  as $j\tendi$.
  We now claim that we can 
  \ben
\item change the normalisation from $\|\gamma_D^- w^-_j\|_{L^2(\Gamma)}=1$ to $\|\gamma_D^- w^-_j\|_{H^{1/2}_{|k_j|}(\Gamma)}=|k_j|^{1/2}$ (or indeed any finite power of $|k_j|$), and 
\item assume, without loss of generality, that $k_j \in \Rea$ for all $j$.
  \een
  Indeed, \cite[Corollary 6.1]{MoSp:19} shows that the results of \cite{PoVo:99} imply existence of a quasimode normalised by 
  $\|\gamma_D^- w^-_j\|_{H^{1}_{|k_j|}(\Gamma)}=|k_j|$, and then \cite[Corollary 6.2]{MoSp:19} shows that this implies existence of a quasimode  with $k_j\in \Rea$ for all $j$, normalised by  $\|\gamma_D^- w^-_j\|_{H^{1}_{k_j}(\Gamma)}=k_j$.
To obtain the claim involving Points 1 and 2 above, we need to justify that we can replace the normalisation 
  $\|\gamma_D^- w^-_j\|_{H^{1}_{|k_j|}(\Gamma)}=|k_j|$ by $\|\gamma_D^- w^-_j\|_{H^{1/2}_{|k_j|}(\Gamma)}=|k_j|^{1/2}$. This follows by repeating the arguments in \cite[Corollaries 6.1 and 6.2]{MoSp:19} with
  \bit
\item the bound on the Dirichlet-to-Neumann map from $H_k^1(\Gamma)\rightarrow L^2(\Gamma)$ from \cite[Lemma 5]{MoLu:68} replaced by the analogous bound from $H^{1/2}_k(\Gamma)\rightarrow H^{-1/2}_k(\Gamma)$ obtained by interpolation (see, e.g., \cite[Lemma 4.2]{ChSpGiSm:20}), and 
\item the bounds on the $L^2(\Gamma)\rightarrow L^2_{\rm comp}(\Rea^d$) norms of $\cV_{i/o}$ and $\cK_{i/o}$ from \cite[Lemma 4.3]{Sp2013a} replaced by analogous bounds on the $H^{-1/2}_k(\Gamma)\rightarrow L^2_{\rm comp}(\Rea^d)$ and 
  $H^{1/2}_k(\Gamma)\rightarrow L^2_{\rm comp}(\Rea^d)$ norms, respectively; these bounds are proved using the bounds in Lemma \ref{lem:Newton}, the trace result \eqref{eq:Dtrace}, and similar arguments to those in the proof of Lemma \ref{lem:Neumann}.
  \eit
Note that, in both these points, the precise algebraic powers of $k_j$ don't matter, since they are dominated by the $O(k_j^{-\infty})$ coming from the quasimode.

  With the changes to the quasimode $w_j^\pm$ in Points 1 and 2 above, we now let
  \begin{align}\label{eq:QMv-}
    v_j^- := w_j^- + \cN_i \big((\Delta + k^2_j n_{i}) w_j^{-}\big) \quad\tand\quad 
    v_j^+ := w_j^+ + \cN_o \big((\Delta + k^2_j n_{o}) w_j^{+}\big) 
  \end{align}
  (where the arguments of $\cN_i$ and $\cN_o$ are extended by zero outside their supports),
  and observe that, since  $w_j^+$ has compact support, $v_j^+ \in \SRC(k_j \sqrt{n_o})$. 
  Let 
  \begin{equation*}
    {\bf g}_j := \Cauchy^- v_j^- -  \Cauchy^+ v_j^+;
  \end{equation*}
  then $v$ satisfies \eqref{eq:BVP2} with $c_i=n_i$ and $c_o=n_o$ and ${\bf f= g}_j$.

  We now show that \eqref{eq:QM0} holds. On the one hand, by the definition of ${\bf g}_j$, the second equation in \eqref{eq:QM1}, Lemma \ref{lem:Newton}, and the first equation in \eqref{eq:QM1}.
  \begin{align}\nonumber
    \N{{\bf g}_j}_{H^{1/2}_{k_j}(\Gamma)\times H^{-1/2}_{k_j}(\Gamma)} 
    &\leq \N{\Cauchy^- w_j^- -  \Cauchy^+ w_j^+}_{H^{1/2}_{k_j}(\Gamma)\times H^{-1/2}_{k_j}(\Gamma)} 
    + \N{\Cauchy^- \cN_i \big((\Delta + k^2_j n_{i}) w_j^{-}\big)}_{H^{1/2}_{k_j}(\Gamma)\times H^{-1/2}_{k_j}(\Gamma)} \\ \nonumber
    &\quad+ \N{\Cauchy^+ \cN_o \big((\Delta + k^2_j n_{o}) w_j^{+}\big)}_{H^{1/2}_{k_j}(\Gamma)\times H^{-1/2}_{k_j}(\Gamma)},\\
    &= O(k_j^{-\infty}) \quad\tas j\tendi. \label{eq:comb1}
  \end{align}
  On the other hand, using \eqref{eq:QMv-}, the normalisation $\|\gamma_D^- w_j^-\|_{H^{1/2}_{k_j}(\Gamma)}=k_j^{1/2}$, 
  Lemma \ref{lem:Newton}, and \eqref{eq:QM1} (in that order), we have
  \begin{align}\nonumber
    \N{\Cauchy^- v_j^-}_{H^{1/2}_{k_j}(\Gamma)\times H^{-1/2}_{k_j}(\Gamma)} \geq \N{\gamma_D^- v_j^-}_{H^{1/2}_{k_j}(\Gamma)} 
    & \geq \N{\gamma_D^- w_j^-}_{H^{1/2}_{k_j}(\Gamma)} - \N{\gamma_D^- \cN_i \big((\Delta + k^2_j n_{i}) w_j^{-}\big)}_{H^{1/2}_{k_j}(\Gamma)},\\
    & = k^{1/2}_j + O(k_j^{-\infty}) \quad\tas j\tendi.\label{eq:comb2}
  \end{align}
  The bound \eqref{eq:QM0} (and hence the result) follows by combining \eqref{eq:comb1} and \eqref{eq:comb2}.
  \epf

  \bpf[Proof of Theorem \ref{thm:bound2}]
  This follows by combining Theorem \ref{thm:1}, Lemma \ref{lem:bound2}, and Lemma \ref{lem:bound3}.
  Indeed, if $n_i<n_o$, then $S_{io}$ has ``good'' behaviour via Lemma \ref{lem:bound2}, but $S_{oi}$ has ``bad'' behaviour via Lemma \ref{lem:bound3}. 
  If $n_i>n_o$, then $S_{io}$ has ``bad'' behaviour via Lemma \ref{lem:bound3}, and $S_{oi}$ has ``good'' behaviour via Lemma \ref{lem:bound2}. 
  \epf

  We also record the following upper bound on $\N{\First^{-1}}_{H^{1/2}_k\times H^{-1/2}_k}$, valid for all Lipschitz $\Oi$.
  
  %
  %
  \begin{theorem}\mythmname{Inverse is algebraically bounded in frequency for almost all frequencies}\label{thm:bound1}
    Given positive real numbers $n_i,n_o, k_0, \delta,$ and $\varepsilon$, there exists a set $J\subset [k_0,\infty)$ with $|J|\leq \delta$ 
    and $C=C(\delta,\varepsilon,k_0)$ 
    such that
    \begin{equation*}
      \N{\First^{-1}}_{H^{1/2}_k\times H^{-1/2}_k} \leq C k^{2 + 5d/2+\varepsilon} \quad\tfa k\in [k_0,\infty)\setminus J.
    \end{equation*}
  \end{theorem}
  
  This result shows that, at high frequency, the blow-up associated with spurious quasi-resonances is extremely localised in frequency, thus giving a possible reason why spurious quasi-resonances seem to have rarely been noticed in literature.


  \bpf[Proof of Theorem \ref{thm:bound1}]
  This follows from Theorem \ref{thm:1}, Corollary \ref{cor:bound1}, and the results of \cite{LaSpWu:20}. 
  Indeed, \cite[Theorem 1.1]{LaSpWu:20} implies that, for arbitrary positive real numbers $n_i,n_o$, the assumptions of Lemma \ref{lem:bound1} (and hence also Corollary \ref{cor:bound1}) are satisfied with  $C_{\rm sol}(k) \sim k^{5d/2+1+\varepsilon}$. To see this, we note that \cite[Theorem 1.1]{LaSpWu:20} holds for problems fitting in the ``black-box scattering'' framework, and the transmission problem fits in this framework by \cite[Lemma 2.3 and Remark 2.4]{LaSpWu:20}. Furthermore, \cite[Theorem 1.1]{LaSpWu:20} is an $L^2\rightarrow L^2$ bound, but this implies a bound of the form \eqref{eq:Csol} with $C_{\rm sol}(k) \sim k^{5d/2+1+\varepsilon}$ thanks to Green's identity -- see the comments around \cite[Equation 1.3]{LaSpWu:20}.
  \epf


  \section{Proofs of Lemma \ref{lem:aug1} and Theorem \ref{thm:infsup} (the  results about the augmented BIEs)}\label{sec:proof_aug}

  \bpf[Proof of Lemma \ref{lem:aug1}]
  We first prove the result when $\widetilde{A}_* =\AugFirst$.
  By the first equation in \eqref{eq:augmented} $\bphi = \First^{-1}{\bf g}$. Then, by the second equation in \eqref{eq:augmented} and the expression for $\First^{-1}$ \eqref{eq:T1S1}, 
  \begin{equation*}
    {\bf 0} = P_i^+ \bphi = P_i^+ \First^{-1} {\bf g} = P_i^+ (S_{io}+S_{oi}-I) {\bf g}.
  \end{equation*}
  By \eqref{eq:L41} $P_i^+ S_{io}=0$, and then 
  \begin{equation}\label{eq:constraint2}
    P^+_i(S_{oi}-I){\bf g} = {\bf 0}.
  \end{equation}
  Now, by \eqref{eq:L42}, $P^-_i(S_{oi}-I){\bf g} = {\bf 0}$, and thus the constraint \eqref{eq:constraint} follows by \eqref{eq:sum}.
  Then, by 
  the first equation in \eqref{eq:augmented} and \eqref{eq:T1S1}, 
  $\bphi=\First^{-1}{\bf g} = (S_{io} + S_{oi}-I){\bf g}$, 
  and the result \eqref{eq:augsol} follows by \eqref{eq:constraint}.

  We now prove the result when $\widetilde{A}_* =\AugSecond$. Similar to before, by \eqref{eq:augmented} and the expression for $\Second^{-1}$ \eqref{eq:T1S2},
  \begin{equation*}
    {\bf 0} = P_i^+ \bphi = P_i^+ \Second^{-1} {\bf g} = P_i^+ (I-S_{io}-S_{oi}+2 S_{io}S_{oi}) {\bf g}= P^+_i(I-S_{oi}){\bf g}, 
  \end{equation*}
  since $P_i^+ S_{io}=0$ by \eqref{eq:L41}; we therefore obtain \eqref{eq:constraint2}, and the constraint \eqref{eq:constraint} follows exactly as before.
  The result \eqref{eq:augsol} then follows from using the constraint \eqref{eq:constraint} in the expression for $\Second^{-1}$ \eqref{eq:T1S2}.
  \epf

  \bpf[Proof of Theorem \ref{thm:infsup}]
  We first prove \eqref{eq:infsup1}.
  Let $\bpsi= (\bpsi_1,\bpsi_2)$ for $\bpsi_1, \bpsi_2\in \cH$. Then
  \begin{equation*}
    \big(\AugFirst \bphi,\bpsi\big)_{\cH\times\cH} = \big(\First \bphi, \bpsi_1\big)_{\cH} + \big( P_i^+ \bphi,\bpsi_2\big)_{\cH}.
  \end{equation*}
  Given $\bphi$, let $\bpsi_2:= P_i^+ \bphi$. Since $\bphi= P_i^+\bphi+ P_i^-\bphi$,
  \begin{equation*}
    \big(\AugFirst \bphi,\bpsi\big)_{\cH\times\cH} = \big(\First P_i^+\bphi, \bpsi_1\big)_{\cH} +
    \big(\First P_i^-\bphi, \bpsi_1\big)_{\cH}
    + \N{P_i^+ \bphi}^2_{\cH}.
  \end{equation*}
  Motivated by Theorem \ref{thm:physical}, let $\bpsi_1:= S_{io}^* P_i^-\bphi$. 
  By Theorem \ref{thm:physical}, $S_{io}\First =I$ as an operator $R(P_i^-)\rightarrow R(P_i^-)$, and thus
  \begin{align*}
    \big(\AugFirst \bphi,\bpsi\big)_{\cH\times\cH} &
    = \big(\First P_i^+\bphi,\, S_{io}^* P_i^-\bphi\big)_{\cH} +
    \big(S_{io}\First P_i^-\bphi,P_i^-\bphi\big)_{\cH}
    + \N{P_i^+ \bphi}^2_{\cH}\\
    &= \big(\First P_i^+\bphi,\, S_{io}^* P_i^-\bphi\big)_{\cH} +
    \N{P_i^-\bphi}_{\cH}^2
    + \N{P_i^+ \bphi}^2_{\cH}.
  \end{align*}
  Now, by the second equality in \eqref{eq:First}, Lemma \ref{lem:0}, and \eqref{eq:sum},
  \begin{align*}
    \big(\First P_i^+\bphi, \,S_{io}^* P_i^-\bphi\big)_{\cH} 
    &= \big(S_{io}(P_i^- -P_o^+)P_i^+\bphi, \,P_i^-\bphi\big)_{\cH},
  \end{align*}
  which equals zero since $P_i^- P_i^+=0$ by \eqref{eq:sum} and $S_{io}P_o^+=0$ by Lemma \ref{lem:SP}.
  Therefore, with this choice of $\bpsi$, 
  \begin{align}
    \frac{
      \big|
      \big( \AugFirst \bphi, \bpsi\big)_{\cH\times\cH}
      \big|
    }{
      \N{\bphi}_{\cH} \N{\bpsi}_{\cH\times\cH}
    }
    &= \frac{
      \N{P_i^- \bphi}^2_{\cH}+\N{P_i^+ \bphi}^2_{\cH}
    }{
      \N{\bphi}_{\cH} \sqrt{ \N{S_{io}^* P_i^- \bphi}^2_{\cH} + \N{P_i^+\bphi}^2_{\cH}}
    }
    \geq 
    \frac{
      \sqrt{\N{P_i^- \bphi}^2_{\cH}+\N{P_i^+ \bphi}^2_{\cH}}
    }{
      \N{\bphi}_{\cH} \max\big\{ \N{S_{io}}_{\cH\rightarrow\cH}, 1 \big\}
    }
    .\label{eq:infsup2}
  \end{align}
  By \eqref{eq:sum}, the triangle inequality, and the inequality $2ab\leq a^2 + b^2$ for $a,b>0$, 
  \begin{equation}\label{eq:infsup3}
    \N{P_i^- \bphi}^2_{\cH}+\N{P_i^+ \bphi}^2_{\cH}\geq \frac{1}{2} \N{\bphi}^2_{\cH}.
  \end{equation}
  The result \eqref{eq:infsup1} then follows from combining \eqref{eq:infsup2} and \eqref{eq:infsup3}.

  We now prove \eqref{eq:infsup1a}.
  As above, let $\bpsi= (\bpsi_1,\bpsi_2)$ for $\bpsi_1, \bpsi_2\in \cH$. Motivated by the proof of \eqref{eq:infsup1}, given $\bphi \in \cH$, let $\bpsi_1:= S_{io}^* P_i^- \bphi$. Then, by  the definition of $\AugSecond$ and Theorem \ref{thm:physical},
  \begin{align}\nonumber
    \big(\AugSecond \bphi,\bpsi\big)_{\cH\times\cH} &
    = \big(S_{io}\Second \bphi,\, P_i^-\bphi\big)_{\cH} + \big(P_i^+\bphi, \bpsi_2\big)_{\cH},\\
    &= \big(S_{io}\Second P_i^+\bphi,\, P_i^-\bphi\big)_{\cH} 
    + \N{P_i^- \bphi}^2_{\cH}+
    \big(P_i^+\bphi, \bpsi_2\big)_{\cH}.
    \label{eq:infsup4}
  \end{align}
  Now, by the definition of $\Second$ \eqref{eq:Second}, Lemma \ref{lem:0}, \eqref{eq:sum}, and the fact that  $S_{io}P_o^+=0$ by Lemma \ref{lem:SP},
  \begin{equation}\label{eq:infsup5}
    \big(S_{io} \Second P_i^+\bphi,\, P_i^- \bphi\big)_{\cH} = 
    \big(S_{io}\big( P_o^- + I\big) P_i^+\bphi, \, P_i^-\bphi\big)_{\cH} =
    2\big(S_{io} P_i^+\bphi, \, P_i^-\bphi\big)_{\cH}.
  \end{equation}
  We now let $\bpsi_2:= P_i^+\bphi - 2 S_{io}^* P_i^- \bphi$. This definition along with \eqref{eq:infsup4} and \eqref{eq:infsup5} implies that
  \begin{align*}
    \big( \AugSecond \bphi, \psi\big)_{\cH\times \cH} &= 2 \big(S_{io} P_i^+ \bphi,\, P_i^- \bphi\big)_{\cH} + \N{P_i^- \bphi}^2_{\cH} 
    + \N{P_i^+ \bphi}^2_{\cH}
    - 2 \big( P_i^+ \bphi, S_{io}^* P_i^- \bphi\big)_{\cH},\\
    &= \N{P_i^- \bphi}^2_\cH+ \N{P_i^+ \bphi}^2_\cH.
  \end{align*}
  Therefore, with this choice of $\bpsi$,
  \begin{align}\nonumber
    \frac{
      \big|
      \big( \AugSecond \bphi, \bpsi\big)_{\cH\times\cH}
      \big|
    }{
      \N{\bphi}_{\cH} \N{\bpsi}_{\cH\times\cH}
    }
    = \frac{
      \N{P_i^- \bphi}^2_{\cH}+\N{P_i^+ \bphi}^2_{\cH}
    }{
      \N{\bphi}_{\cH} \sqrt{ \N{S_{io}^* P_i^- \bphi}^2_{\cH} + \N{(P_i^+-2 S_{io}^* P_i^-)\bphi}^2_{\cH}}
    }.
  \end{align}
We now use the triangle inequality and \eqref{eq:Cauchy} to find that 
\begin{align*}
&\N{S_{io}^* P_i^- \bphi}^2 + \N{(P_i^+-2 S_{io}^* P_i^-)\bphi}^2
\leq  \N{S_{io}^*}^2\N{ P_i^- \bphi}^2 + \Big(\N{P_i^+\bphi} + 2\N{S_{io}^*}\N{ P_i^-\bphi}\Big)^2\\
&\hspace{3cm}\leq  \N{S_{io}^*}^2\N{ P_i^- \bphi}^2+ \N{P_i^+\bphi}^2 + 4\N{S_{io}^*}^2\N{ P_i^-\bphi}^2 + 4  \N{P_i^+\bphi}
\N{S_{io}^* }\N{P_i^-\bphi}\\
&\hspace{3cm}\leq (1+ 2\epsilon) \N{P_i^+\bphi}^2+(5 + 2/\epsilon) \N{S_{io}^*}^2\N{ P_i^- \bphi}^2_{\cH};
\end{align*}
if $\epsilon= 1+ \sqrt{2}$, then $5+2\epsilon=1+2\epsilon = 3+2 \sqrt{2}$ and thus
\beqs
    \frac{
      \big|
      \big( \AugSecond \bphi, \bpsi\big)_{\cH\times\cH}
      \big|
    }{
      \N{\bphi}_{\cH} \N{\bpsi}_{\cH\times\cH}
    }
       \geq 
    \frac{
      \sqrt{\N{P_i^- \bphi}^2_{\cH}+\N{P_i^+ \bphi}^2_{\cH}}
    }{
      \N{\bphi}_{\cH}\sqrt{3+ 2\sqrt{2}} \max\big\{ \N{S_{io}}_{\cH\rightarrow\cH}, 1 \big\}
    }.
\eeqs
The result \eqref{eq:infsup1a} follows from this inequality and \eqref{eq:infsup3}.
  \epf
  
  \appendix
  
  \section{Extension of the results to the more general form of the transmission problem}\label{app:general}

We now sketch how the decomposition \eqref{eq:T1S1} of the first-kind BIE extends to the more general transmission problem of Definition \ref{def:HTP} with the transmission condition in \eqref{eq:BVP} replaced by \eqref{eq:newtc}. The analogous extension of the decomposiiton \eqref{eq:T1S2} for the second-kind BIE is very similar.

\paragraph{Derivation of the first kind BIE.}
The equation \eqref{eq:L31} holds as before, but now the analogue of \eqref{eq:L34} is 
\beqs
P_o^- D^{-1} \Cauchy^- u^- = P_o^- D^{-1} {\bf f}.
\eeqs
Therefore, the analogue of the first-kind BIE $\First$ in \eqref{eq:firstsecondkind} is
\beqs
\Firstgen \Cauchy^- u^- = P_o^- D^{-1}{\bf f}
\eeqs
where
\begin{align}\label{eq:Firstgen}
  \Firstgen
  &:= P^-_oD^{-1}-D^{-1}P^+_i = D^{-1}P^-_i - P^+_oD^{-1}=
  \left[
    \begin{array}{cc}
      -(K_i+ K_o) & V_i+ \alpha^{-1}V_o\\
      \alpha^{-1}W_i+ W_o &\alpha^{-1}( K'_i+ K'_o)
    \end{array}
  \right].
\end{align}

\paragraph{Solution operators.}
Let $S(c_i,c_o)$ be the solution operator of the BVP \eqref{eq:BVP2} with the transmission condition replaced by \eqref{eq:newtc}. 
Let $\widetilde{S}(c_i,c_o)$ be the solution operator to \eqref{eq:BVP2} with the transmission condition 
\beqs
\Cauchy^- u^- = D^{-1} \Cauchy^+ u^+  +{\bf f};
\eeqs
i.e., $\alpha$ is replaced by $\alpha^{-1}$ compared to $S(c_i,c_o)$. Let $S_{io}:= S(n_i,n_o)$ and let $\widetilde{S}_{oi}:= \widetilde{S}(n_o,n_i)$.

\ble
\beq\label{eq:genformula}
\big( \Firstgen\big)^{-1} = S_{io} D + D \widetilde{S}_{oi} -D.
\eeq
\ele

\bpf[Sketch of the proof]
The analogue of Lemma \ref{lem:4} is now that 
  \begin{equation}\label{eq:L41gen}
    \bphi = S_{io}{\bf f} \quad\text{ if and only if } 
    \quad
    \left\{
      \begin{array}{c}
        P_i^- \bphi=\bphi, \tand\\
        P_o^-D^{-1}(\bphi -{\bf f}) ={\bf 0}.
      \end{array}
    \right.
  \end{equation}
and
  \begin{equation}\label{eq:L42gen}
    \bphi = \widetilde{S}_{oi}{\bf f} \quad\text{ if and only if } 
    \quad
    \left\{
      \begin{array}{c}
        P_o^- \bphi=\bphi, \tand\\
        P_i^-D(\bphi -{\bf f}) ={\bf 0}.
      \end{array}
    \right.
  \end{equation}
We then repeat the steps in the proof of Theorem \ref{thm:1}. 
Assuming that $ \Firstgen\bpsi={\bf g}$ and applying $P_i^-D$, 
we find that the analogue of \eqref{eq:T11a} is that
\beqs
P_i^-D\big(P_o^-D^{-1} \bpsi - {\bf g}\big) = {\bf 0},
\eeqs
so that, by \eqref{eq:L42gen},
\beq\label{eq:end2}
P_o^- D^{-1} \bpsi = \widetilde{S}_{oi} {\bf g}.
\eeq
Similarly applying $P_o^-$ to $\Firstgen\bpsi={\bf g}$, we find that the analogue of \eqref{eq:T12a} is 
\beqs
P_o^-\big(D^{-1}P_i^-\bpsi - {\bf g}\big) = {\bf 0},
\eeqs
so that, by \eqref{eq:L41gen},
\beq\label{eq:end3}
P_i^- \bpsi = S_{io} D{\bf g}.
\eeq
Then, using that $ \Firstgen\bpsi={\bf g}$, \eqref{eq:Firstgen}, \eqref{eq:end2}, and \eqref{eq:end3}, we find that the analogue of \eqref{eq:end1} is 
\beqs
\bpsi = (P_i^+ + P_i^-)\bpsi = D P_o^- D^{-1}\psi - D {\bf g} + P_i^- \bpsi = D \widetilde{S}_{oi} {\bf g} - D{\bf g} +S_{io} D{\bf g},
\eeqs
which is \eqref{eq:genformula}.
\epf

  \section*{Acknowledgements}

  EAS thanks Zo\"is Moitier (KIT) for telling him about the references \cite{CaWi:15} and \cite{LeDjArLaZyDuScBo:07}, and Alastair Spence (University of Bath) for useful discussions about augmented operator equations. 
    We thank the referees for their careful reading of the paper and constructive comments.
  EAS was supported by EPSRC grant  EP/R005591/1. AM acknowledges support from GNCS-INDAM, from PRIN project ``NAFROM-PDEs'' and from MIUR through the ``Dipartimenti di Eccellenza'' Program (2018-2022) - Dept.~of Mathematics, University of Pavia.

  \bibliographystyle{plain}

\end{document}